\documentclass[psamsfonts]{amsart}
\usepackage{amssymb, amscd, amsmath}
\usepackage{diagrams}
\diagramstyle[scriptlabels,small]

\DeclareMathOperator{\adj}{adj}
\DeclareMathOperator{\Der}{Der}
\DeclareMathOperator{\id}{id}
\DeclareMathOperator{\tr}{tr}
\DeclareMathOperator{\cok}{cok}
\DeclareMathOperator{\Hom}{Hom}
\DeclareMathOperator{\Ext}{Ext}
\DeclareMathOperator{\End}{End}

\DeclareMathOperator{\Spec}{Spec}

\DeclareMathOperator{\depth}{depth}

\DeclareMathOperator{\grade}{grade}
\DeclareMathOperator{\ann}{Ann}
\DeclareMathOperator{\pd}{pd}
\DeclareMathOperator{\rk}{rank}
\DeclareMathOperator{\rank}{\rk}
\DeclareMathOperator{\Cl}{Cl}
\DeclareMathOperator{\diag}{diag}
\DeclareMathOperator{\syz}{syz{}}
\DeclareMathOperator{\Alt}{Alt{}}
\DeclareMathOperator{\ev}{ev}

\newcommand{\p}{{\mathfrak{p}}}

\newcommand{\m}{{\mathfrak{m}}}

\renewcommand{\phi}{\varphi}
\renewcommand{\bar}{\overline}
\renewcommand{\to}{\longrightarrow}
\newcommand{\Cal}{\mathcal}
\newcommand{\Sc}{{\Cal Sc}}
\newcommand{\J}{{\Cal J}}

\renewcommand{\tilde}{\widetilde}
\renewcommand{\P}{P}
\newcommand{\uHom}{\underline{\Hom}}
\newcommand{\uEnd}{\underline{\End}}
\newcommand{\uExt}{\underline{\Ext}}
\newcommand{\calE}{\mathcal{E{}}}

\newcommand{\bbbd}{{\mathbb D}}
\newcommand{\bbbh}{{\mathbb H}}

\newcommand{\bbbn}{{\mathbb N}}

\newcommand{\bbbs}{{\mathbb S}}

\newcommand{\bbbz}{{\mathbb Z}}

\newcommand{\oF}{{\overline {F}}}
\newcommand{\oG}{{\overline {G}}}

\newcommand{\oN}{{\overline {N}}}

\newarrow{Equal}{}{=}{}{=}{}
\newarrow{Dotted}{}{.}{}{.}{>}
\newarrow{Dashed}{}{dash}{}{dash}{>}

\def\xto{\xrightarrow}

\theoremstyle{plain}
\newtheorem{theorem}{Theorem}
\numberwithin{theorem}{section}

\newtheorem{prop}[theorem]{Proposition}

\newtheorem{lemma}[theorem]{Lemma}

\newtheorem{cor}[theorem]{Corollary}

\newtheorem*{matrixmainthm}{Theorem~\ref{thm:matrixmain}}
\newtheorem*{n_is_3thm}{Theorem~\ref{thm:n_is_3}}
\newtheorem*{wild}{Corollary~\ref{cor:wild}}
\newtheorem*{normfactclassify}{Theorem~\ref{thm:normfactclassify}}
\newtheorem*{corfactorizations}{Corollary~\ref{cor:factorizations}}

\theoremstyle{definition}
\newtheorem{defn}[theorem]{Definition}

\newtheorem{remark}[theorem]{Remark}
\newtheorem{sit}[theorem]{}
\newtheorem{notation}[theorem]{Notation}

\newtheorem{example}[theorem]{Example}

\newtheorem{prob}[theorem]{Problem}

\numberwithin{equation}{theorem}

\begin{document}

\title[MCM Modules on the Determinant]{%
Factoring the Adjoint and\\
Maximal Cohen--Macaulay Modules \\
over the Generic Determinant} 

\author[R.-O. Buchweitz]{Ragnar-Olaf Buchweitz}
\address{Dept.\ of Comp. and Math. Sciences, 
University of Toronto at Scarborough (UTSC), 
Toronto, Ont.\ M1A 1C4, Canada}
\email{ragnar@math.utoronto.ca}

\author[G.J. Leuschke]{Graham J. Leuschke}
\address{Dept.\ of Math., Syracuse University,
Syracuse NY 13244, USA}
\email{gjleusch@math.syr.edu}
\urladdr{http://www.leuschke.org/}

\thanks{The first author was partly supported by NSERC grant
3-642-114-80, and the second author by an NSA Young Investigator grant.}

\date{\today}


\begin{abstract}
A question of Bergman~\cite{Bergman:2003} asks whether the adjoint of
the generic square matrix over a field can be factored nontrivially as
a product of square matrices.  We show that such factorizations indeed
exist over any coefficient ring when the matrix has even size.
Establishing a correspondence between such factorizations and
extensions of maximal Cohen--Macaulay modules over
the generic determinant, we exhibit all factorizations where one of
the factors has determinant equal to the generic determinant.  The
classification shows not only that the Cohen--Macaulay representation
theory of the generic determinant is wild in the tame-wild dichotomy,
but that it is quite wild: even in rank two, the isomorphism classes
cannot be parametrized by a finite-dimensional variety over the
coefficients. We further relate the factorization problem to the
multiplicative structure of the $\Ext$--algebra of the two nontrivial
rank-one maximal Cohen--Macaulay modules and determine it completely.

\end{abstract}
\maketitle


\section{Introduction}

Let $K$ be a field, $X = (x_{ij})$ the generic $(n\times n)$--matrix, 
whose entries thus form a family of $n^2$ indeterminates, and $S =
K[x_{ij}]$, the polynomial ring over $K$ in those variables.  The
determinant $\det X$ of $X$ is a homogeneous polynomial of degree $n$
with coefficients $\pm 1$, and the hypersurface
ring $R = S/(\det X)$ is a normal domain of dimension $n^2-1$. 

The classical adjoint, or adjugate, $\adj(X)$ of $X$ is uniquely
determined through either of
the two matrix equations
\begin{equation}
\label{eq:defadj}
\adj(X)X = (\det X)\cdot\id_{n}\quad\text{and}\quad X\adj(X)= 
(\det X)\cdot\id_{n}\,, 
\end{equation}
where $\id_n$ denotes the $(n\times n)$ identity matrix.

G.M.~Bergman asks \cite{Bergman:2003} whether the factorizations
(\ref{eq:defadj}) and those arising from the transposes $X^T$,\,
$\adj(X)^T$ are the 
only nontrivial matrix factorizations of $(\det X)\cdot\id_n$\,.  More
specifically, he  
inquires about possible refinements of the factorization 
(\ref{eq:defadj}) obtained by writing $\adj(X) = YZ$ for noninvertible
$(n\times n)$--matrices $Y$ and $Z$.  
He shows, for $K$ an algebraically closed field of characteristic
zero, that there are no such refinements when $n$ is odd and that for
$n$ even the only possible refinements have either $\det Y = \det X$
or $\det Z = \det X$, up to multiplication by units in $S$.  
The proofs of \cite{Bergman:2003} use a recent theorem by C.~De
Concini and Z.~Reichstein \cite{DeConcini-Reichstein:2003} about maps
between Grassmannians, generalizing the well-known topological theorem
that the hairy sphere cannot be combed.

Here, in Section~\ref{section:matrix}, we show that when $n$ is even,
the adjoint can in fact be factored nontrivially (over any
commutative ring $K$).   We give explicit matrix factorizations for
each invertible \emph{alternating\/} matrix $A$ over $S$, based on the
following key result:

\begin{matrixmainthm}
Let $U,A$ be $(n\times n)$--matrices over a commutative ring $K$ with
$A$ alternating and $\det U$ a nonzerodivisor in $K$.  There exist
then unique alternating $(n\times n)$--matrices $B_A$ and ${}_AB$
satisfying
\[
A \adj(U) = U^T B_A \qquad\text{ and }\qquad \adj(U) A ={}_AB U^T\,.
\]
\end{matrixmainthm}
When $n$ is even, there exist \emph{invertible\/} alternating matrices
$A$, so that $Y = A^{-1}U^T$, $Z = B_A$ gives one of the
factorizations of $\adj(U)$ allowed by Bergman's result.
\begin{corfactorizations}
Let $X$ be the generic square matrix of even size over the commutative
ring $K$. Then the adjoint $\adj(X)$ factors nontrivially. 
\end{corfactorizations}

The remainder of the paper has two main purposes: to show that the
factorizations arising from this Corollary are the only
factorizations possible with $\det Y = \det X$ or $\det Z =\det X$ up
to units, and to cover, in reverse, the path by which we
found them.  To this end, we observe (Proposition~\ref{prop:corresp})
that a factorization of $\adj(X)$ into a product of
square matrices $Y$ and $Z$ exhibits the cokernel of $\adj(X)$ as 
the middle term in a short exact sequence of $R$-modules, with ends
the modules presented by $Z$ and $Y$.  Each of the
three modules in this extension is a maximal Cohen--Macaulay
$R$-module, and so is given by a matrix factorization of $\det X$.  In
Section~\ref{section:matfacts} we briefly discuss the essential
features that we will need from the theory of matrix factorizations. 

Bergman's question can thus be rephrased in terms of extensions: Is it
possible to write the cokernel of the adjoint as an extension of two
maximal Cohen--Macaulay $R$-modules?  When $K$ is a unique
factorization domain, W.~Bruns 
has shown \cite{Bruns:1975} (see also \cite{Bruns-Roemer-Wiebe:2005})
that up to isomorphism there
are only three MCM $R$-modules of rank one, namely the cokernel of
$X$, the cokernel of the transpose $X^T$, and $R$ 
itself.  This observation, together with a calculation in 
the divisor class group of $R$, already allows us to give a negative
answer to the $n=3$ case of Bergman's question over any UFD.

\begin{n_is_3thm}
Let $X = (x_{ij})$ be the generic $(3\times3)$--matrix over a unique
factorization domain $K$.  Then there are no nontrivial factorizations
of $\adj(X)$. 
\end{n_is_3thm}

The general question of identifying whether and under what conditions
a given module can be the middle term of a nonsplit short exact
sequence is interesting and rarely addressed.  We avoid it here as
well.  Looking instead for inspiration to Bergman's theorem we observe
that the condition $\det Y = u\det X$, with $u$ a unit, is equivalent
to the module 
presented by $Y$, $\cok Y$, having rank one as an $R$-module.  Given
the classification of rank-one MCM $R$-modules, we obtain an explicit
correspondence between nontrivial factorizations $\adj(X) = YZ$ with
$\det Y = u\det X$ and short exact sequences 
\[
  \begin{diagram}[midshaft]
  0&\rTo &\cok Y &\rTo&Q&\rTo&\cok X&\rTo& 0
  \end{diagram}
\]
such that $Q$ is a homomorphic image of $R^n$ (see
Lemma~\ref{lem:minexts-surjs} and Proposition~\ref{prop:our-extn}).  

The classification of factorizations $\adj(X)=YZ$ with $\det Y = u\det
X$ thus naturally leads to the calculation of $\Ext^1$ for
the rank-one MCM $R$-modules. In Section~\ref{section:ext1-L-L} we
show that $\Ext_R^1(\cok X, \cok 
X)=0$.  This follows from a theorem of R.~Ile \cite{Ile:2004}; we
reprove Ile's result, simplifying the proof slightly.  We compute the
minimal graded free resolution of $\Ext_R^1(\cok X, \cok X^T)$ in
Theorem~\ref{thm:res-E}. 

Sections~\ref{section:hom-M-Lstar} and \ref{section:L-Lstar} classify
the nontrivial factorizations $\adj(X)=YZ$ with $\det Y = u\det X$ and
the associated extensions.  We show 
\begin{normfactclassify}
Let $\adj(X) = YZ$ be a factorization of $\adj(X)$ with $\det Y = u
\det X$ for some unit $u$.  Then $\cok Y \cong \cok X^T$ and $Y=
JX^TZ$ for an invertible $(n \times n)$--matrix $J$.  Moreover, there
exist then a unique invertible alternating $(n \times n)$--matrix $A$
and a matrix $U$ of the same size such that 
\[
J^{-1} = A + X^TU {\text{ \qquad and \qquad }} Z = B_A + U\adj(X)\,.
\]
Two such factorizations $\adj(X)=JX^TZ$ and $\adj(X)=J'X^TZ'$ give the
same extension if and only if $J^{-1} - J'^{-1} = X^TV$ for some
$(n\times n)$--matrix $V$, and in that case $Z-Z' = V \adj(X)$.  
\end{normfactclassify}
As explained above, this result depends upon the structure of
$\Ext_R^1(\cok X, \cok X^T)$ determined in
Section~\ref{section:L-Lstar}. 

Given the classification of the maximal Cohen--Macaulay
$R$-modules of rank $1$, one may ask for a description of the maximal
Cohen--Macaulay modules of small rank in general.  From the results of
earlier sections, in Section~\ref{section:extsofrk1s} we make a first
step in this direction by classifying all extensions of the rank-one
maximal Cohen--Macaulay modules.  In representation-theoretic terms,
the class of such extension modules is (very) \emph{wild\/}: 

\begin{wild}
Let $X=(x_{ij})$ be the generic $(n \times n)$--matrix over a
field $K$, $n\geq 3$.  Let $R = K[x_{ij}]/(\det X)$ be the
generic determinantal hypersurface ring.  Then the rank-two
orientable MCM $R$-modules cannot be parametrized by the points of any
finite-dimensional algebraic variety over $K$.
\end{wild}

Finally, we construct a graded ring $\calE$, the \emph{stable
$\Ext$-algebra\/} of the rank-one maximal Cohen--Macaulay modules and
describe its multiplication, given by the Yoneda product, explicitly.
This algebra controls the higher-order extension theory of the 
rank-one MCM $R$-modules.  

Our thanks go to George Bergman for several interesting email
exchanges about this material and related ideas.  We are also pleased
to acknowledge our debt to the computer algebra system {\tt
  Macaulay2} \cite{M2}.

\section{The Adjoint of even size factors}\label{section:matrix}

In this section we give the promised factorizations of $\adj(X)$,
after some background on determinants and derivations.  Throughout,
$K$ denotes a commutative ring, $X = (x_{ij})$ the generic $(n\times
n)$--matrix over $K$, and $S = K[x_{ij}]$\,. 

\begin{sit}
We will use the following notation for minors of the generic
matrix $X$: Let $[i_{1}i_{2}\cdots i_{k}\mid
j_{1}j_{2}\cdots j_{k}]$ denote the (unsigned) determinant
of the $(k\times k)$--submatrix of $X$ that consists of the
rows indexed $1\leq i_{1}<\cdots < i_{k}\leq n$, and of the
columns indexed $1\leq j_{1}<\cdots < j_{k}\leq n$.

The symbol $[i_{1}i_{2}\cdots i_{k}\;\widehat\mid\;
j_{1}j_{2}\cdots j_{k}]$ will denote the complementary
minor, thus, the determinant of the
$(n-k)\times(n-k)$--submatrix of $X$ obtained by removing
the rows indexed $i_{\nu}$ and the columns indexed
$j_{\nu}$.  For consistency, the empty determinant $[\
  \;\mid\;\ ]$, for
$k=0$, has value $1$, whereas the empty complementary minor $[\
  \;\widehat\mid\;\ ]$ equals $\det X$.

We extend the symbols $[?\mid\; ?]$ and $[?\;\widehat\mid\;?]$
to not necessarily strictly increasing index sets by
requiring them to be alternating in both the left and right
arguments. In particular, each symbol vanishes if there is 
repetition of indices either before or after the vertical 
bar.  
\end{sit}

\begin{sit}\label{sit:eval}
If $U$ is any $(n\times n)$--matrix over some
$K$-algebra $R$, then there exists a unique $K$-algebra
homomorphism $\ev_{U}:S\to R, x_{ij}\mapsto u_{ij}$, that
transforms the entries of $X$ to those of $U$.  The evaluation
homomorphism $\ev_U$ is compatible with the formation of minors, thus
$[\cdots](U) = \ev_{U}([\cdots])$ represents the corresponding 
minor of the matrix $U$.  Of particular interest is the case $U=X^T$.
The corresponding evaluation homomorphism $\tau:= \ev_{X^T}$ is then a
$K$-algebra involution of $S$, given by $\tau(x_{ij}) = x_{ji}$, that
fixes the determinant as $\det \tau X = \tau(\det X) = \det X$.

We write $I_{t}(U)\subseteq R$ for
the ideal generated by all the $(t\times t)$--minors of $U$.
The transpose of a matrix $U$ will be denoted $U^{T}$.  We also
sometimes write $|U| := \det U$ to abbreviate.
\end{sit}

\begin{example}
The $(i,j)^{\text{th}}$ entry of the adjoint matrix can be written as 
\begin{equation*}\label{eq:adjentries}
\adj(X)_{ij} =
(-1)^{i+j}[j\;\widehat\mid\; i] = (-1)^{i+j}[1\cdots\widehat
j\cdots n\mid 1\cdots\widehat i\cdots n]\,.
\end{equation*}
\end{example}

\begin{sit}
Recall that a map $D:R\to R$, on a not necessarily
commutative ring $R$, is a {\em derivation\/} if it satisfies the
Leibnitz rule $D(ab) =
D(a)b + aD(b)$ for all elements $a,b\in R$.

For example, the {\em partial derivative\/} $\partial_{ij} =
\frac{\partial}{\partial x_{ij}}$ with respect to the
variable $x_{ij}$ defines a derivation on $S$ that is
furthermore $K$-linear.  These partial derivations form
indeed a basis of the free $S$-module $\Der_{K}(S)$ of all
$K$-linear derivations on $S$,
$$
\Der_{K}(S) \cong \bigoplus_{1\le i,j\le n}S\partial_{ij}\,.
$$
\end{sit}

Now we state the facts on derivations and minors that we 
will use.

\begin{lemma}
\label{lem:diffdet}
If $R$ is a commutative ring, $D:R\to R$ a derivation,
and $U$ an $(n\times n)$--matrix over $R$, then $D(\det U)$ can be
written as a sum of determinants,
\begin{equation*}\label{eq:diffdet}
D(\det U) = \sum_{i=1}^{n}
\left|
\begin{matrix}
u_{11}&\cdots&u_{1n}\\
\vdots&&\vdots\\
D(u_{i1})&\cdots&D(u_{in})\\
\vdots&&\vdots\\
u_{n1}&\cdots&u_{nn}
\end{matrix}
\right|
= \sum_{i=1}^{n}
\left|
\begin{matrix}
u_{11}&\cdots&D(u_{1j})&\cdots&u_{1n}\\
\vdots&&\vdots&&\vdots\\
u_{n1}&\cdots&D(u_{nj})&\cdots&u_{nn}
\end{matrix}
\right|
\,.
\end{equation*}
\end{lemma}
\begin{proof}
This follows immediately from the Leibnitz rule for
$D$ applied to the complete expansion of the
determinant.
\end{proof}
\begin{lemma}
\label{lem:derdet} 
Let $X$ be again the generic matrix and $S$ the associated
polynomial ring over $K$.
\begin{enumerate}
\item[(1)]
\label{lem:derdet1}
For any pair of indices $1\le i,j\le n$,
$$
\partial_{ij}(\det X) = \adj(X)_{ji}\,,
$$ 
equivalently,
$$
\adj(X)^{T} = (\partial_{ij}(\det X) )_{ij}\,.
$$

\item[(2)]
\label{lem:derdet2}
For any pair of indices $1\le i,j\le n$,
$$
\sum_{\nu=1}^{n}x_{i\nu}\partial_{j\nu}(\det X) =
\delta_{ij}\det X  = \sum_{\nu=1}^{n}x_{\nu i}\partial_{\nu
j}(\det X)\,,
$$
where $\delta_{ij}$ is the {\em Kronecker symbol\/}.
\item[(3)]
\label{lem:derdet3}
For any indices $1\le i_{1}, i_{2},\ldots,i_{k}\le n$ and $1\le
j_{1},j_{2},\ldots,j_{k}\le n$, 
\begin{align*}
\partial_{i_{1}j_{1}}\cdots \partial_{i_{k}j_{k}}(\det X) &=
(-1)^{i_{1}+\cdots +i_{k}+j_{1}+\cdots j_{k}} [i_{1}\cdots
i_{k}\;\widehat\mid\, j_{1},\cdots,j_{k}]
\,;
\end{align*} 
in particular, these terms vanish whenever there is a 
repetition among the $i$'s or the $j$'s.
\end{enumerate}
\end{lemma}

\begin{proof}
Claim (1) follows from Lemma~\ref{lem:diffdet} with
$D=\partial_{ij}$ and $U=X$.  In view of (1), claim (2) is
simply a reformulation of the equation (\ref{eq:defadj}) above.
To see (3), apply first Lemma~\ref{lem:diffdet} or (1) to the generic
matrix using the derivation $\partial_{i_{k}j_{k}}$, and 
then use induction on $k\ge 1$.
\end{proof}

\begin{sit}
We now use the ``differential calculus'' from above to establish two
factorization results about 
products of the adjoint matrix with {\em
alternating\/} matrices on one or both sides.  Recall that
an $(n\times n)$--matrix $A=(a_{kl})$ is alternating if
$A^{T}= -A$ and the diagonal elements vanish, $a_{kk}=0$ for
each $k=1,\ldots, n$.  The latter condition is of course a
consequence of the first as soon as $2$ is a nonzerodivisor
in $K$.
\end{sit}

\begin{theorem}
\label{thm:matrixmain}
Let $U, A$ be $(n\times n)$--matrices over a commutative
ring $K$, with $A$ alternating.  The $(n\times n)$--matrix
$B_A = (b_{rs})$ with entries from $I_{1}(A)\cdot I_{n-2}(U)
\subseteq K$, given by
\begin{equation*}\label{eq:Bentries}
b_{rs}
= \sum_{k<l}a_{kl}(-1)^{r+s+k+l}[rs\;\widehat\mid\; kl](U)\,,
\end{equation*}
is then alternating as well and satisfies the matrix
equation
\begin{equation}\label{eq:Bprop}
A \adj(U) = U^T B_A\,.
\end{equation}
If $\det U$ is a nonzerodivisor in $K$, then $B_A$ is the {\em
unique\/} solution to this equation.
\end{theorem}

\begin{proof}
As $[sr\;\widehat\mid\; kl] = -[rs\;\widehat\mid\; kl]$ and
$[rr\;\widehat\mid\; kl] = 0$, the matrix $B=B_A$ is 
alternating.  To verify that $B$ satisfies (\ref{eq:Bprop}), it
suffices to establish the generic case, in which we replace $K$ by $S$ and
$U$ by $X$.  Let $E_{ij}$ denote the {\em elementary\/} $(n\times
n)$--matrix with $1$ at position $(i,j)$ as its only nonzero
entry.  Recall that $E_{ab}E_{cd}= \delta_{bc}E_{ad}$ for
any indices $1\le a,b,c,d \le n$.  As
$\partial_{rk}\partial_{sl}(\det X) =
(-1)^{r+s+k+l}[rs\;\widehat\mid\; kl]$ by
Lemma~\ref{lem:derdet}(3), the right-hand side of
(\ref{eq:Bprop}) expands first as
\begin{align*}
X^TB &= \bigg(\sum_{i,\nu}x_{\nu i}E_{i\nu}\bigg)
\bigg(\sum_{\mu,j}\sum_{k<l}a_{kl}\partial_{\mu k}
\partial_{jl}(\det X)E_{\mu j}\bigg)\\
&=\sum_{k<l}a_{kl}\sum_{i,j}\bigg(\sum_{\nu}x_{\nu i}
\partial_{\nu k}\partial_{jl}(\det X) \bigg)E_{ij}\,.
\end{align*}
The innermost sum can be simplified using first that partial
derivatives commute, then applying the product
rule, and finally invoking Lemma~\ref{lem:derdet}(2)
together with the fact that $\partial_{jl}(x_{\nu i}) =
\delta_{j\nu}\delta_{li}$. In detail, these steps yield the 
following equalities:
\begin{align*}
\sum_\nu x_{\nu i} \partial_{\nu k}\partial_{jl}(\det X)
&= \sum_\nu x_{\nu i} \partial_{jl}\partial_{\nu k}(\det X)\\
&= \sum_\nu \partial_{jl}\big(x_{\nu i} \partial_{\nu k}(\det
X)\big) - \sum_\nu
\partial_{jl}(x_{\nu i})\partial_{\nu k}(\det X) \\
&= \partial_{jl}\bigg(\sum_\nu x_{\nu i} \partial_{\nu k}(\det
X)\bigg) - \delta_{li} \sum_\nu \delta_{j\nu}
\partial_{\nu k}(\det X) \\
&= \delta_{ik} \partial_{jl} (\det X) -
  \delta_{li} \partial_{jk}(\det X)\,.
\end{align*}
In light of this simplification, we may expand $X^TB$ further 
as follows:
\begin{align*}
X^TB &=\sum_{k<l}a_{kl}\sum_{i,j}\bigg(\sum_{\nu}x_{\nu i}
\partial_{\nu k}\partial_{jl}(\det X) \bigg)E_{ij}\\
&= \sum_{k<l}a_{kl}\sum_{i,j}\bigg(\delta_{ik}
\partial_{jl}(\det X) - \delta_{li}\partial_{jk}(\det
X)\bigg)E_{ij}\\
&= \sum_{k<l}a_{kl}\sum_{j}\bigg(\partial_{jl}(\det X)E_{kj}
- \partial_{jk}(\det X)E_{lj}\bigg)\\
&= \sum_{k<l}a_{kl}\bigg(E_{kl}\sum_{j}\partial_{jl}(\det
X)E_{lj} -E_{lk} \sum_{j} \partial_{jk}(\det
X)E_{kj}\bigg)\\
&= \sum_{k<l}a_{kl}\bigg(E_{kl}\sum_{i,j}
\partial_{ji}(\det X)E_{ij} -E_{lk}
\sum_{i,j}\partial_{ji}(\det X)E_{ij}\bigg)\\
&= \sum_{k<l}a_{kl}\big(E_{kl}-E_{lk}\big)
\sum_{i,j}\partial_{ji}(\det X)E_{ij}\\
&= A \adj(X)
\end{align*}
with the last equality using that $A$ is alternating, thus
$A=\sum_{k<l}a_{kl}(E_{kl}-E_{lk})$, and that $\adj(X) =
\sum_{i,j}\partial_{ji}(\det X)E_{ij}$, in view of
Lemma~\ref{lem:derdet}(1).

The final assertion about uniqueness follows from
(\ref{eq:Bprop}) by multiplying from the left with $\adj(U)^T$ and
using equation (\ref{eq:defadj}) to obtain
\[
\adj(U)^T A \adj(U) = (\det U) \cdot B\,.
\]
\end{proof}

\begin{cor}\label{cor:otherside}
For $U$, $A$ as in Theorem~\ref{thm:matrixmain}, there exists also an
alternating $(n \times n)$--matrix ${}_AB$ so that 
\[
\adj(U)A = {}_AB U^T\,.
\]
\end{cor}

\begin{proof}
Let $B_A$ be the $(n \times n)$--matrix over $S = K[x_{ij}]$ given by
Theorem~\ref{thm:matrixmain}, so that $A\adj(X) = X^TB_A$,  and let
$\tau = \ev_{X^T}$ be the  
involution introduced in \ref{sit:eval}.  Clearly $\tau$ exchanges $X$
and its transpose, and moreover, $\tau(\adj(X)) = \adj(X)^T$, in view
of equation (\ref{eq:defadj}).  Now
\begin{alignat*}{2}
A\adj(X) &= X^TB_A&\quad&\text{if, and only if,}\\
\tau(A)\tau(\adj(X)) &=\tau(X^T)\tau(B_A)&&\text{if, and only if,}\\
\tau(A)\adj(X)^T &= X\tau(B_A)&&\text{if, and only if,}\\
\adj(X)\tau(A)^{T} &= \tau(B_A)^{T}X^T\,.
\end{alignat*}
As $A$ and $B_A$ are both alternating, so are $\tau(A)$ and
$\tau(B_A)$, and the last equation is equivalent to
\[
\adj(X)\tau(A) = \tau(B_A)X^T\,.
\]
Interchanging the roles of $A$ and $\tau(A)$, we have 
\[
\adj(X) A = \tau(B_{\tau(A)}) X^T\,.
\]
Put ${}_AB = \tau(B_{\tau(A)}$.
\end{proof}

We now investigate what happens when multiplying simultaneously from
both left and right. 

\begin{prop}
\label{prop:matrixmain2}
Let $U,A$ be again $(n\times n)$--matrices over a commutative ring
$K$, and let $B_A$ be the matrix introduced in
Theorem~\ref{thm:matrixmain}.  For another alternating $(n\times 
n)$--matrix $A' = (a'_{uv})$, the $(n\times n)$--matrix $C=C_{A,A'}=(c_{wm})$
with entries from $I_{1}(A)\cdot I_{n-3}(U)\cdot I_{1}(A')\subseteq R$ 
given by
\begin{align*}
c_{wm}&= \sum_{k<l, u<v}(-1)^{u+v+w+k+l+m}a_{kl}
[uvw\;\widehat\mid\;klm](U)a'_{uv}
\end{align*}
{satisfies}
\begin{equation*}
B_A A' = r\cdot\id_{n} + CU^T \qquad\text{ and }\qquad A{}_{A'}B =
r\cdot \id_n + U^T C\,,
\end{equation*}
{where}
\begin{align*}
r &= -\sum_{k<l, u<v}(-1)^{u+v+k+l}a_{kl}
[uv\;\widehat\mid\;kl](U)a'_{uv}\ \in K\,.
\end{align*}
\end{prop}

\begin{proof}
It suffices again to verify the result for the generic 
matrix $U=X$, in which case we can employ once more the 
description of minors as given in Lemma~\ref{lem:derdet}(3).
The straightforward calculation proceeds then as follows:
\begin{align*}
(B_A A' - r\cdot\id_{n})_{ij}&=
\sum_{m}\sum_{k<l}(-1)^{i+m+k+l}a_{kl}[im\;\widehat\mid\;
kl]a'_{m j} \\
&\quad\quad + \delta_{ij}\sum_{k<l}\sum_{u<v}(-1)^{u+v+k+l}a_{kl}
[uv\;\widehat\mid\;kl]a'_{uv}\\
&=\sum_{k<l}a_{kl}\left(\sum_{m}\partial_{ik} 
\partial_{ml}(\det X)a'_{m j} + \sum_{u<v}\partial_{uk} 
\partial_{vl}(\det X)a'_{uv}\delta_{ij}\right)\\
&=\sum_{k<l,m}a_{kl}\left(\partial_{ik}
\partial_{ml}(\det X)a'_{mj} + \sum_{u<v}\partial_{uk}
\partial_{vl}\left(\partial_{im}(\det X)x_{jm}\right)a'_{uv}\right)
\end{align*}
where we have used \ref{lem:derdet}(2) in the last step.
Using the product rule twice together with
$
\partial_{ab}(x_{cd})=\delta_{ac}\delta_{bd}\,,
$
we find next
\begin{align*}
\partial_{uk} \partial_{vl}\left(\partial_{im}(\det X)x_{jm}\right)
&= \partial_{vl} \partial_{im} (\det X) \delta_{uj} \delta_{km} +
\partial_{uk} \partial_{im}(\det X) \delta_{vj}\delta_{lm} \\
& \quad \quad + \partial_{uk} \partial_{vl} \partial_{im}(\det X) x_{jm}\,.
\end{align*}
Substituting and evaluating the Kronecker symbols yields
\begin{align*}
(B_A A' - r\cdot\id_{n})_{ij}
&=\sum_{k<l,m}a_{kl}\left(\partial_{ik}
\partial_{ml}(\det X)a'_{mj} + \sum_{u<v}\partial_{uk}
\partial_{vl}\left(\partial_{im}(\det X)x_{jm}\right)a'_{uv}\right) \\
&=\sum_{k<l}a_{kl}\left(\sum_{m}\partial_{ik} \partial_{ml}(\det X)a'_{m j} +
\sum_{j<v}\partial_{vl} \partial_{ik}(\det X)a'_{jv}\right.\\
&\qquad\quad  \left. + \sum_{u<j}\partial_{uk}
\partial_{il}(\det X)a'_{uj} + \sum_{u<v,m}\partial_{uk}
\partial_{vl}\partial_{im}(\det X)a'_{uv}x_{jm}
\right)
\end{align*}
The terms involving only second order derivatives
of the determinant cancel. To see this, rename
summation indices, and use that
$
\partial_{mk} \partial_{il}(\det
X)=-\partial_{ml} \partial_{ik}(\det X)
$ 
and that $A'$ is alternating, whence its entries satisfy
$a'_{mm}=0, a'_{jm}=-a'_{mj}$.  In detail,
\begin{align*}
(B_A A' - r\id_{n})_{ij}
&=\sum_{k<l}a_{kl}\left(\sum_{m}\partial_{ml}
\partial_{ik}(\det X)a'_{jm} + 
\sum_{j<m}\partial_{ml}
\partial_{ik}(\det X)a'_{jm} \right.\\
&\quad\quad  \left. +
\sum_{m<j}\partial_{mk}
\partial_{il}(\det X)a'_{mj}
+ \sum_{u<v}\sum_{m}\partial_{uk}
\partial_{vl}\partial_{im}(\det X)a'_{uv}x_{jm}
\right)\\
&=\sum_{m}\left(\sum_{k<l,u<v} a_{kl} \partial_{uk}
\partial_{vl}\partial_{im}(\det X)a'_{uv}\right)x_{jm}\\
&=\sum_{m}\left(\sum_{k<l, u<v} (-1)^{u+v+i+k+l+m}a_{kl}
[uvi\;\widehat\mid\;klm]a'_{uv}\right)x_{jm}\\
&= (CX^T)_{ij}
\end{align*}
where we have evaluated the third order derivatives of the
determinant according to \ref{lem:derdet}(3).

For the other statement, we observe that
\begin{align*}
A' {}_AB X^T &= A'\adj(X)A\\
&= X^T B_{A'}A\\
&= X^T (r\cdot\id_n + CX^T)\,,
\end{align*}
so that 
\[
A' {}_AB = r\cdot\id_n + X^T C
\]
as well.
\end{proof}

Combining the results from Theorem~\ref{thm:matrixmain}
and Proposition~\ref{prop:matrixmain2} yields the following.

\begin{theorem}
\label{thm:bifac}
Let $U,A,A'$ be $(n\times n)$--matrices over a commutative
ring $K$, with $A,A'$ {\em alternating\/}. One then has an 
equality of matrices
$$
A\adj(U)A' = r\cdot U^T + U^T C U^T\,,
$$
where $r$ and $C=C_{A,A'}$ are as specified in {\em 
Proposition~\ref{prop:matrixmain2}\/}.  In particular, $A\adj(U) A'$ is
both left- and right-divisible by $U^T$. 
\qed
\end{theorem}

\begin{remark}\label{rmk:half-trace}
The element $r\in I_{1}(A)\cdot I_{n-2}(U)\cdot I_{1}(A')\subseteq R$
is 
a ``{\em half trace\/}'' of both $B_A A'$ and $A {}_{A'}B$, as
\begin{align*}
\tr(B_A A') &= \sum_{k<l}\sum_{i,j}a_{kl}(-1)^{i+j+k+l}[ij\;
\widehat\mid\;kl]a'_{ji}\\
&= 2\sum_{k<l}\sum_{i<j}a_{kl}(-1)^{i+j+k+l}[ij\;
\widehat\mid\;kl]a'_{ji}\\
&=2r
\end{align*}
invoking once again that $A'$ is alternating.  Equivalently, 
$\tr(CU^T) = (2-n)r\,.$
\end{remark}
\medskip

\begin{remark}\label{rmk:n=2}
If $n=2$, all expressions of the form
$[uvw\;\widehat\mid\;klm]$ vanish, so that Theorem~\ref{thm:bifac}
(for $U=X$) specializes to the easily established identity
\begin{align*}
\begin{pmatrix}
0&a\\
-a&0
\end{pmatrix}
\begin{pmatrix}
x_{22}&-x_{12}\\
-x_{21}&x_{11}
\end{pmatrix}
\begin{pmatrix}
0&b\\
-b&0
\end{pmatrix}
=
-ab
\begin{pmatrix}
x_{11}&x_{21}\\
x_{12}&x_{22}
\end{pmatrix}\,.
\end{align*}
\end{remark}
\bigskip

If $n=2m$ is even, then over any commutative ring there are {\em
invertible\/} alternating matrices of size $n$.  For example, the
alternating ``hyperbolic matrix'' 
$
\begin{pmatrix}
0&\id_{m}\\
-\id_{m}&0
\end{pmatrix}
$ has determinant equal to 1 over any ring.  The following corollary,
immediate from Theorem~\ref{thm:bifac}, thus gives the
factorizations promised in the Introduction.

\begin{cor}\label{cor:factorizations}
If $n$ is {\em even\/}, then the adjoint of the 
generic matrix admits nontrivial factorizations
$$
\adj(X) = YZ = Y'Z'
$$
into products of $(n\times n)$--matrices over 
$S$ with $\det Y = \det Z' = \det X$ up to units of $S$.

More precisely, any pair of invertible alternating $(n\times
n)$--matrices $A,A'$ over $S$
gives rise to such factorizations.  With $r$ and $C$ the
data associated to $A, A'$ as in {\em
  Proposition~\ref{prop:matrixmain2}\/}, 
one may take
\begin{align*}
Y &= A^{-1}X^{T}\quad\text{and}\quad Z = 
(r\cdot \id_{n} + CX^{T}){A'}^{-1}\,,\\
Y' &= A^{-1}(r\cdot\id_{n} + X^{T}C)\quad\text {and} \quad Z' =
X^{T}{A'}^{-1}\,. 
\end{align*}
\end{cor}

\section{Matrix Factorizations}\label{section:matfacts}

The rest of this paper is devoted to interpreting the factorizations
of Section~\ref{section:matrix} as extensions of maximal
Cohen--Macaulay modules over the hypersurface ring $R = S/(\det
X)$.  Here we collect some preliminary material, including a brief
r\'esum\'e of the theory of matrix factorizations, after D.~Eisenbud
\cite{Eisenbud:1980}, and some convenient results on stable
homomorphism modules and multilinear algebra.  

\begin{defn}\label{sit:matfactdef}
Let $S$ be a commutative Noetherian ring.  A \emph{matrix
  factorization\/} $(\phi, \psi, F, G)$ of an 
element $f \in S$ is a pair of homomorphisms between finitely
generated free $S$-modules, $\phi: G \to F$ and $\psi: F \to G$,
satisfying $\phi\psi=f\cdot\id_F$ and $\psi\phi = f\cdot\id_G$.  We
sometimes suppress $F$ and $G$ from the notation and refer to the
matrix factorization $(\phi,\psi)$.
\end{defn}

\begin{sit}\label{sit:matfactmods}
Let $(\phi,\psi,F,G)$ be a matrix factorization of $f\in S$, and
 assume that $f$ is a nonzerodivisor.  Then we have exact sequences
  \begin{equation}\label{eq:matfact-cokers}
  \begin{diagram}[midshaft]
  0&\rTo &G&\rTo^{\phi}&F&\rTo&\cok\phi&\rTo& 0\\
  0&\rTo &F&\rTo^{\psi}&G&\rTo&\cok \psi&\rTo& 0\,.
  \end{diagram}
  \end{equation}
As $f \cdot F = \phi\psi(F)$ is contained in the image of $\phi$, the
cokernel of $\phi$ 
is annihilated by $f$.  Similarly, $f \cdot \cok \psi =0$.  Thus $\cok
\phi$ and $\cok \psi$ are naturally finitely generated modules over $R =
S/(f)$.  If we write  $\overline{?} := ? \otimes_S R$ for reduction
modulo $f$, then the sequence 
\begin{equation}\label{eq:matfactres}
\begin{diagram}[midshaft]
\cdots & \rTo^{\bar\psi} & \oG & \rTo^{\bar\phi} & \oF &
\rTo^{\bar\psi} & \oG & \rTo^{\bar\phi} & \oF\ 
& (\rTo & \cok\bar\phi & \rTo & 0\ ) 
 \end{diagram}
\end{equation}
is a complex of free $R$-modules that
constitutes a free resolution of $\cok \phi = \cok \bar\phi$.  In
particular, $\cok \phi$ has a periodic resolution of period at most $2$.

The reversed pair $(\psi,
\phi)$ is also a matrix factorization of $f$.  Put $M = \cok \phi$ and 
$N = \cok \psi$; then the resolution (\ref{eq:matfactres}) 
exhibits $N$ as a first syzygy of $M$ over $R$ and vice versa:
  \begin{equation*}
  \begin{diagram}[midshaft]
  0&\rTo &N&\rTo&\oF&\rTo&M&\rTo& 0\\
  0&\rTo &M&\rTo&\oG&\rTo&N&\rTo& 0
  \end{diagram}
  \end{equation*}
are exact sequences of $R$-modules.
As a matter of notation, we write $M = \cok(\phi,\psi)$ and $N =
\cok(\psi,\phi)$ to emphasize their provenance.  There are two
distinguished \emph{trivial\/} matrix factorizations, 
namely $(1,f,S,S)$ and $(f,1,S,S)$.  Note that $\cok(1,f)=0$, while
$\cok(f,1) \cong R$.   
\end{sit}

\begin{sit}\label{sit:ranks}
Suppose again that $f \in S$ is a nonzerodivisor.  Then the free
modules $F$ and $G$ in any matrix factorization $(\phi,\psi,F,G)$ of
$f$ have the same rank $n$, as can be seen from
equation~(\ref{eq:matfact-cokers}).  The homomorphisms $\phi$ and
$\psi$, then, can be represented by square matrices over $S$ after
choosing bases for $F$ and $G$.

If in addition $f$ is a prime element of $S$, so that $R$ is an
integral domain, then from
$\phi\psi=f\cdot \id_n$ it follows that both $\det \phi$ and $\det \psi$
are, up to units, powers of $f$.  Specifically, $\det \phi = u f^k$
and $\det \psi = u^{-1}f^{n-k}$ for some unit $u \in S$ and $k\leq
n$.  In this case the $R$-module $\cok (\phi,\psi)$ has
rank $k$, while $\cok(\psi,\phi)$ has rank $n-k$.  (To see this,
localize at the prime ideal $(f)$.  Then over the discrete valuation
ring $S_{(f)}$, $\phi$ is equivalent to $f\cdot \id_k \oplus
\id_{n-k}$ and so $\cok\phi$ has rank $k$ over the field $R_{(f)}$.)
\end{sit}

\begin{defn}\label{sit:matfacthoms}
Given two matrix factorizations $(\phi_1,\psi_1,F_1,G_1)$ and
$(\phi_2,\psi_2,F_2,G_2)$ of the same element $f\in S$, a homomorphism
of matrix   
factorizations from $(\phi_1,\psi_1)$ to $(\phi_2,\psi_2)$ is a pair
of homomorphisms of free modules $\alpha: F_1 \to F_2$ and $\beta: G_1
\to G_2$ rendering 
\begin{equation*}
\begin{diagram}[midshaft]
F_1 & \rTo^{\psi_1} & G_1 & \rTo^{\phi_1} & F_1 \\
\dTo<{\alpha} & & \dTo>{\beta} & & \dTo>{\alpha} \\
F_2 & \rTo^{\psi_2} & G_2 & \rTo^{\phi_2} & F_2 
\end{diagram}
\end{equation*}
commutative.  Such a diagram induces a homomorphism of $R$-modules
$\cok(\phi_1,\psi_1) \to \cok(\phi_2,\psi_2),$ which we write as
$\cok(\alpha,\beta)$. 

Consider the pullback square
\begin{equation*}
\begin{diagram}[midshaft,width=4em]
C & \rTo & \Hom_S(F_1,F_2) \\
\dTo & \qquad\qquad\SEpbk\  & \dTo>{?\phi_1}  \\
\Hom_S(G_1,G_2) & \rTo_{\phi_2?} & \Hom_S(G_1,F_2)\,.
\end{diagram}
\end{equation*}
The module $\Hom_S((\phi_1,\psi_1),(\phi_2,\psi_2)) := C$ consists of pairs 
\[
(\alpha,\beta) \in \Hom_S(F_1,F_2) \times \Hom_S(G_1,G_2)
\]
so that
$\alpha\phi_1=\phi_2\beta$.  There is a natural map $\Hom_S(F_1,G_2)
\to C$ sending $\gamma \in
\Hom_S(F_1,G_2)$ to 
$(\phi_2\gamma,\gamma\phi_1)$, and an exact sequence
\begin{equation*}
\begin{diagram}[midshaft]
\Hom_S(F_1,G_2) & \rTo & C & \rTo & \Hom_R(M_1,M_2) & \rTo & 0\,,
\end{diagram}
\end{equation*}
which is also exact at the left if $\phi_2$ is injective,
\emph{e.g.\/} if $f$ is a nonzerodivisor in $S$.

The two matrix factorizations are {\em equivalent\/} if there is a
homomorphism of matrix factorizations $(\alpha,\beta)$ as above in
which both 
$\alpha$ and $\beta$ are isomorphisms of free modules.  Direct sums of 
matrix 
factorizations are defined in the natural way, and we say that a
matrix factorization $(\phi,\psi)$ is {\em reduced\/} provided it is
not equivalent to a matrix factorization with a direct
summand of the form $(f,1)$.  Similarly, $(\phi,\psi)$
is called \emph{minimal\/} if it is not equivalent to one with a
direct summand of the form $(1,f)$.
\end{defn}

The following Proposition, while straightforward to verify, is key in
our constructions of matrix factorizations.

\begin{prop}\label{sit:pushout}
Let $(\alpha,\beta): (\phi_1,\psi_1,F_1,G_1) \to
(\phi_2,\psi_2,F_2,G_2)$ be a homomorphism of matrix factorizations of 
$f \in S$, set $R = S/(f)$, 
and put $M_i =  
\cok(\phi_i,\psi_i), N_i = \cok(\psi_i,\phi_i)$ for $i=1,2$.  Then the 
bottom row of the pushout diagram of $R$-modules
\begin{equation}\label{eq:pushout}
\begin{diagram}[midshaft]
0 & \rTo &  M_1 &  \rTo & \overline{G_1} & \rTo & N_1 & \rTo & 0 \\
  & & \dTo<{\cok(\alpha,\beta)}   & \NWpbk\qquad & \dTo & & \dEqual \\
0 & \rTo & M_2 & \rTo & Q & \rTo & N_1 & \rTo & 0
\end{diagram}
\end{equation}
defines an element of $\Ext_R^1(N_1,M_2),$ which is the image of
$\cok(\alpha,\beta)$ under the natural surjection $\Hom_R(M_1,M_2)
\to \Ext_R^1(N_1,M_2).$  The module $Q$ is again 
given by a matrix factorization, namely
$$
Q \cong \cok\left(
\left(\begin{matrix} \phi_2 & \alpha \\ 0 & \psi_1\end{matrix}\right), 
\left(\begin{matrix} \psi_2 & -\beta \\ 0 & \phi_1\end{matrix}\right)
\right)\,.
$$
\end{prop}

\begin{sit}\label{sit:zero-Ext}
If, in the notation of \ref{sit:pushout}, $\cok(\alpha,\beta)$ factors
through a projective $R$-module, then the bottom row of
(\ref{eq:pushout}) splits, and vice versa.  In this case,
$\cok(\alpha,\beta)$ factors through $\oG_1$, and we have $Q \cong M_2
\oplus N_1$. 
\end{sit}

The main application we have in mind for matrix factorizations is
their equivalence with maximal Cohen--Macaulay modules over a
hypersurface ring.  

\begin{defn}\label{def:MCMdef}
Let $R$ be a Noetherian ring and $M$ a finitely generated $R$-module.
Recall that $M$ is a \emph{Cohen--Macaulay\/} module provided
$\depth_{R_\p}M_\p = \dim M_\p$ for each prime $\p \in \Spec R$.  In
particular, $M$ is \emph{maximal Cohen--Macaulay\/} if $M$ is
Cohen--Macaulay and $\dim M = \dim R$.
\end{defn}

\begin{sit}\label{sit:matfactMCMs}
To describe the connection between matrix factorizations and MCM
modules, we let $S$ be a regular ring over which all projective
modules are free.  ({In general, one must replace $F$ and $G$ by
$S$-projectives, and $\phi, \psi$ by appropriate linear maps.})  Let
$f \in S$ be a nonzero nonunit and set $R = S/(f)$.  Given a matrix
factorization $(\phi,\psi,F,G)$ of $f$, we have seen that
$\cok(\phi,\psi)$ has projective dimension $1$ over $S$.  By the Depth
Lemma, we obtain
\[
\depth_R\cok(\phi,\psi)_\p= \dim R_\p
\]
for each $\p \in \Spec R$, so that $\cok(\phi,\psi)$ is a MCM
$R$-module. 

Conversely, let $M$ be a nonzero MCM $R$-module.  Then $\pd_S M = 1$,
so that $M$ has a projective resolution of the form 
\begin{equation}\label{eq:freeresM}
\begin{diagram}[midshaft]
0 & \rTo & G & \rTo^{\phi} & F & \rTo & M & \rTo & 0\,,
\end{diagram}
\end{equation}
with $G$ and $F$ free $S$-modules of the same finite rank.  As $M$ is
annihilated by $f$, the map of complexes from (\ref{eq:freeresM}) to
itself given by multiplication by $f$ is homotopic to zero.
Equivalently, there is a homomorphism $\psi: F \to G$ so that
$\phi\psi=f\cdot \id_F$.  Since $\phi$ is necessarily injective, we have
$\psi\phi=f\cdot \id_G$ as well.  Thus $(\phi,\psi,F,G)$ is a matrix
factorization of $f$ with $\cok(\phi,\psi) \cong M$.
\end{sit}

The matrix factorization $(\phi,\psi,F,G)$ is reduced if and only if
the $R$-modules $M$ and $N$ are \emph{stable\/}, that is, have no nonzero free
direct summand.  Equivalently, no entry of $\phi$ or $\psi$ is a unit.

\begin{theorem}[{\cite[Theorem~6.3]{Eisenbud:1980}}]\label{thm:eisenbud}
Let $S$ be a regular ring such that all projective $S$-modules are
free and set $R=S/(f)$ for a nonzero nonunit $f$.  The association
\[
(\phi,\psi,F,G) \longleftrightarrow \cok(\phi,\psi)
\]
induces a bijection between equivalence classes of reduced matrix
factorizations of $f$ and isomorphism classes of stable MCM $R$-modules.
\end{theorem}

Among other things, this theorem implies that MCM modules over
the ring $R$ above have periodic resolutions of period at most $2$.
In particular, the modules $\Ext_R^i(M,N)$, for $M$ a MCM
module, are periodic in $i$.  To make this
notion more precise, as well as for later use, we recall the
definition of ``stable homomorphisms''. 

\begin{defn}\label{defn:stable-hom}
Let $M$ and $N$ be finitely generated modules over a ring
$R$.  Denote by $\P(M,N)$ the set of $R$-homomorphisms from $M$ to $N$
that factor through a projective $R$-module, and put 
\[
\uHom_R(M,N) = \Hom_R(M,N)/\P(M,N)\,.
\]
We call $\uHom_R(M,N)$ the \emph{stable $\Hom$--module\/}.
We also write $\uEnd_R(M)$
for $\uHom_R(M,M)$, and refer to it as the \emph{stable endomorphism
  ring\/}. 
\end{defn}

Note that $\P(M,N)$ is the image of the natural homomorphism 
\[
q: N \otimes_R \Hom_R(M,R) \to \Hom_R(M,N)
\]
defined by $q(y\otimes f)(x) = y\cdot f(x)$ for  $y \in
N$, $f \in \Hom_R(M,R)$, and $x \in M$.

In order to have a uniform notation for the periodicity of $\Ext$ over
hypersurface rings, we introduce the \emph{ad hoc\/} notion of stable
extension groups.  See also \cite{Avramov-Buchweitz}. 

\begin{defn}\label{defn:stable-Ext}
Let $M$ and $N$ again be finitely generated modules over a ring $R$.
Define the \emph{stable $\Ext$ groups\/} of $M$ by $N$ by
\[
\uExt_R^i(M,N) = \begin{cases}
                         \uHom(M,N) & \text{ if $i=0$} \\ 
                         \Ext_R^i(M,N) & \text{ if $i>0$.}
                 \end{cases}
\]
\end{defn}

For MCM modules over a hypersurface ring, it follows from
\ref{sit:zero-Ext} and the structure of projective resolutions that
the stable $\Ext$ groups are periodic:

\begin{prop}\label{prop:stable-ext-periodic}
Let $M$ and $N$ be finitely generated modules over a hypersurface ring
$R$, with $M$ MCM.  Then 
$\uExt_R^i(M,N) \cong \uExt_R^{i+2}(M,N)$ for all $i\geq 0$. 
\end{prop}

To facilitate explicit computations, we review some conventions from
the dictionary translating matrices to multilinear algebra --- and
back.  In the statements and proofs to follow, we will use these two
languages interchangeably; while the latter is perhaps more elegant,
the former makes for faster and more transparent calculations.

\begin{sit} \label{sit:basic-matrices}
Let $S$ be a ring.  If $F,G$ are finite free modules with, say, $\rank
F = n$ and $\rank G=m$, then an element $a \in F\otimes_SG^*$ may be
viewed as an $(n\times m)$--matrix over $S$.  Namely, in terms of
given ordered bases $(f_1, \dots, f_n)$ for $F$ and $(g_1, \dots,
g_m)$ for $G$, the element $a$ can be written $a = \sum_{i,j} a_{ij}f_i
\otimes g^*_j$, where $(g_1^*, \dots, g_m^*)$ is the canonical dual
basis for $G^* = \Hom_S(G,S)$, and $A=(a_{ij})$ represents the desired
matrix.  As a matrix, $A$ gives a homomorphism $A: G \to F$.
Equivalently, one may view $a$ as a linear form $\alpha: G\otimes_SF^*
\to S$, with $a_{ij} = \alpha(g_j \otimes f_i^*)$.
\end{sit}

\begin{sit}\label{sit:alt-symm-triang}
In the same vein, an element of the second exterior power $\Lambda^2F$
can be identified with an alternating $(n\times n)$--matrix, that is, an
element of $\Alt_n(S)$.  Recall that a square matrix $A = (a_{ij})$ is
\emph{alternating\/} provided $A^T=-A$ and the diagonal elements
vanish, $a_{ii}=0$.  The canonical projection $F\otimes_SF \to
\Lambda^2F$ becomes in terms of matrices the map $A \mapsto A-A^T$.
The kernel of this epimorphism is again a free $S$-module, denoted
$\bbbd_2 F$.  Its elements can be viewed as the \emph{symmetric\/}
matrices, $A=A^T$.  Continuing with this point of view, the canonical
inclusion $\Lambda^2F \to F\otimes_SF$ views an alternating matrix $A$
simply as a matrix, and the cokernel of this map, denoted $\bbbs_2 F$,
can be identified with the module of all $(n\times n)$--matrices
modulo the alternating ones, or equivalently with the free module of
all (at choice: upper or lower) triangular matrices over $S$.
\end{sit}

\begin{sit}
Note that the map $?+?^{T}:F\otimes_S F\to F\otimes_S F$ that sends a
matrix $A$ to $A+A^{T}$ kills all alternating matrices and returns a
symmetric matrix, thus, induces a canonical map of free $S$-modules
$\bbbs_{2}F\to \bbbd_{2}F$ of equal rank ${\binom{n+1}{2}}$.  However,
this map is not an isomorphism if $2$ is not a unit in $S$; rather, it
assigns to the (upper) triangular matrix $U$ the symmetric matrix
$U+U^{T}$ with diagonal entries in the ideal generated by $2$ in $S$.
Thus, the kernel of that map is $(\ann_{S}2)^{n}$ and its cokernel is
$(S/2S)^{n}$. 
\end{sit}

\section{Factorizations and Extensions}\label{section:factn-extn}

In this section we construct an explicit correspondence between
factorizations $\adj(X) = YZ$ and extensions of maximal
Cohen--Macaulay modules over the generic determinantal hypersurface
ring.

\begin{notation}\label{notation} 
Here is the notation for our default situation throughout the rest of
the paper.  Let $K$ be a field, $n$ a positive integer, $X = (x_{ij})$
the generic $(n\times n)$--matrix over $K$, and $S = K[x_{ij}]$.
({Virtually everything below remains true if we assume only that $K$
is a ``Quillen--Suslin regular'' ring, that is, a regular ring over
which all projective modules are free.  To streamline our exposition,
we leave this extension to the interested reader.})  Put $F = S^n$,
the free module of rank $n$, with canonical ordered basis $(f_1,
\dots, f_n)$, and $G = S^n(-1)$ the free $S$-module of the same rank,
but with ordered basis $(g_1, \dots, g_n)$ whose elements are in
degree $1$ with respect to the natural $\bbbn$-grading on $S$.  ({If
an $S$-module is naturally graded, we shall keep track of its grading.
However, not all modules we consider will be graded.})

We write $R := S/(\det X)$, and
$\oN := N\otimes_S R$ for the reduction of an $S$-module $N$ modulo
the determinant.  The hypersurface ring $R$ is a domain of dimension
$n^2-1$, and 
the singular locus of $R$ (in any characteristic) is defined by the
partial derivatives of $\det X$, which by Lemma~\ref{lem:derdet} are
precisely the entries of $\adj(X)$.  The ideal generated by these
entries, $I_{n-1}(X)$, is prime of height $4$ in $S$
\cite[2.5]{Bruns-Vetter:1988}, so that the singular locus
$V(I_{n-1}(X))$ has codimension $3$ in $\Spec R$. In particular, $R$
is regular in codimension one, and so is a normal domain.

For $M$ a Cohen--Macaulay $S$-module of codepth $t$, we set $M^\vee :=
\Ext_S^t(M,S(-n))$.  Note that if $M$ is a MCM (so free) $S$-module,
then $M^\vee \cong \Hom_S(M,S)$, while for a MCM $R$-module $M$, we have
$M^\vee \cong \Hom_R(M,R)$, up to shifts in grading.

We define the $R$-modules $L,M$, respectively $L^\vee,M^\vee$, through the
exact sequences of $S$-modules
  \begin{equation}\label{eq:modulesS}
  \begin{diagram}[midshaft]
  0&\rTo &G&\rTo^{X}&F&\rTo&L&\rTo& 0\\
  0&\rTo &F(-n)&\rTo^{\adj(X)}&G&\rTo&M&\rTo& 0\\
  0&\rTo &F^\vee&\rTo^{X^{T}}&G^\vee&\rTo&L^\vee&\rTo& 0\\
  0&\rTo &G^\vee(-n)&\rTo^{\adj(X)^{T}}&F^\vee&\rTo&M^\vee&\rTo& 0\,;
  \end{diagram}
  \end{equation}
equivalently, one has exact sequences of $R$-modules
  \begin{equation}\label{eq:modulesR}
  \begin{diagram}[midshaft]
  0&\rTo& M &\rTo&\oF&\rTo& L&\rTo& 0\\
  0&\rTo& L(-n) &\rTo&\oG&\rTo& M&\rTo& 0\\
  0&\rTo &M^{\vee}&\rTo&\oG^{\vee}&\rTo&L^\vee&\rTo& 0\\
  0&\rTo& L^{\vee}(-n) &\rTo&\oF^{\vee}&\rTo& M^\vee&\rTo& 0\,.
  \end{diagram}
  \end{equation}
Each of $L, M, L^\vee, M^\vee$ is a MCM $R$-module, with associated matrix
factorizations $(X,\adj(X),F,G)$, $(\adj(X), X, G, F)$, and so on.   By
\ref{sit:ranks}, $L$ and 
$L^\vee$ have rank one over $R$, while $M$ and $M^\vee$ have rank
$n-1$.  Fixing any $n-1$ columns of $X$, the module $L$ is isomorphic
to the ideal generated by the maximal minors of those rows, while
$L^\vee$ is obtained similarly by fixing any $n-1$ columns
\cite[Thm. A2.14]{Eisenbud:book}.  
In particular, $L$ and $L^\vee$ are indecomposable nonfree
$R$-modules.  To see that $M$ and $M^\vee$ are indecomposable as well,
localize and use the fact that a syzygy of an
indecomposable MCM module over a Gorenstein local ring is again
indecomposable \cite[Lemma 1.3]{Herzog:1978}. 

As in \ref{sit:eval}, let $\tau:R \to R$ be the
$K$-algebra involution induced by $\tau(x_{ij}) = x_{ji}$. 
Then $\tau$ induces an autoequivalence on the category of $R$-modules,
which we denote $\tau^*$,  satisfying $\tau^* L \cong L^\vee$ and
$\tau^*M\cong M^\vee$. 
\end{notation}

Here is the basic link between factorizations of the adjoint and MCM
modules. 

\begin{prop}\label{prop:corresp}
Let $Y$ and $Z$ be square matrices over $S$ so that $\adj(X)=YZ$.
Then $\cok Y$ and $\cok Z$ are MCM $R$-modules, and there is a short
exact sequence
\begin{equation}\label{eq:factn-extn}
\begin{diagram}
0 & \rTo & \cok Z & \rTo & M & \rTo & \cok Y & \rTo & 0\,.
\end{diagram}
\end{equation}
Furthermore, $Y$ and $Z$ are noninvertible if and only if the exact
sequence is nonsplit.  In this case, $\cok Y$ and $\cok Z$ are
nonfree $R$-modules of rank at most $n-2$.
\end{prop}

\begin{proof}
We have $\det(\adj(X)) = (\det X)^{n-1}$, and $\det X$ is an
irreducible element of $S$.  It follows that, up to unit factors, both
$\det Y$ and $\det Z$ are powers of $\det X$.  In particular, both $Y$
and $Z$ are one-to-one as linear maps.  From (\ref{eq:defadj}), we
have $YZX = (\det X)\cdot \id_n$, and multiplying on the right by $Y$
gives $YZXY=(\det X)\cdot Y$.  Cancelling $Y$ from the left, we have
$ZXY = (\det X) 
\cdot \id_n$.  Since also $XYZ= (\det X)\cdot \id_n$, the pair
$(Z,XY)$ is a matrix factorization of $\det X$.  Similarly, $(Y, ZX)$
is as well a matrix factorization of $\det X$.  Thus $\cok Y$ and
$\cok Z$ are MCM $R$-modules, whose ranks sum to $n-1$ by
\ref{sit:ranks}. 

Taking the canonical basis for $S^n$, we view the matrix $Y$ as an
$S$-linear homomorphism $S^n \to G$ and $Z$ as a homomorphism $F
\to S^n$.  Thus we have the commutative diagram  
\begin{equation}
\begin{diagram}[midshaft]
0 & \rTo & F & \rTo^{\adj(X)} & G & \rTo & M & \rTo & 0 \\
  &      & \dTo<{Z} &       & \dEqual & & \dDashed \\
0 & \rTo & S^n & \rTo_{Y} & G & \rTo & \cok Y & \rTo & 0
\end{diagram}
\end{equation}
with exact rows.  This induces the homomorphism $M \to \cok Y$, which
is surjective with kernel isomorphic to $\cok Z$ by the Snake Lemma,
giving the  exact sequence (\ref{eq:factn-extn}).

Since $M$ is indecomposable, the sequence (\ref{eq:factn-extn}) splits
if and only if either $\cok Y$ or $\cok Z$ is zero, equivalently, one
of $Y$ and $Z$ is invertible.  Finally, if 
$\cok Y$ is a nonzero free module, then (\ref{eq:factn-extn}) clearly
splits, and so $\cok Y=0$.  Since $R$ is a Gorenstein ring, free
modules are also injective objects in the subcategory of MCM modules, whence
(\ref{eq:factn-extn}) splits as well if $\cok Z$ is free, and then
$\cok Z=0$.
\end{proof}

\begin{sit} \label{sit:focus-on-Y}
According to Bergman's theorem \cite{Bergman:2003}, we should only
hope to find nontrivial factorizations $\adj(X)=YZ$ satisfying either
$\det Y = u\det X$ or $\det Z = u\det X$ for some unit $u \in S$, at
least in characteristic zero.  Further, we shall from now onward omit
mention of $u$, and tacitly assume the phrase ``up to unit factors in
$S$'' where necessary.  With this in mind, from this point on {\bf we
consider only factorizations of the adjoint in which
${\boldsymbol{\det Y = \det X}}$.}  The case $\det Z = \det X$ can be
recovered by applying the transpose and the automorphism $\tau$ of
\ref{sit:eval}: If $\adj(X)=YZ$ with $\det Z = \det X$, then
$\adj(X)^T=Z^TY^T$, and so $\adj(X) = \tau(Z^T)\tau(Y^T)$ is a
factorization with $\det \tau(Z^T) =\det X$.
\end{sit}

\begin{sit}\label{sit:rank-one}
Since we assume $\det Y=\det X$, the MCM $R$-module $\cok Y$ has rank
one by \ref{sit:ranks}.  It is also reflexive, so isomorphic to a
divisorial ideal of $R$.  The divisor class group of $R$ was computed
by Bruns: $\Cl(R) \cong \bbbz$, generated by the class of
$[L]=-[L^\vee]$ (see \cite{Bruns:1975} or \cite[7.3.5]{BH}).
Furthermore, the symbolic powers $L^{(m)}$ representing elements $m[L]
\in \Cl(R)$ are equal to the usual powers $L^m$, and among these, only
$L$ and $L^\vee$ are MCM modules \cite[9.27]{Bruns-Vetter:1988}.  More
generally, if $K$ is only a normal domain, then $R$ is still normal
and $\Cl(R) \cong \Cl(K) \oplus \bbbz$.  Succinctly: when $\Cl(K)=0$,
the only nonfree MCM $R$-modules of rank one, up to isomorphism, are
$L$ and $L^\vee$.  Thus $\det Y = \det X$ implies either $\cok Y \cong
L$ or $\cok Y \cong L^\vee$.
\end{sit}

This already allows us to rule out all nontrivial factorizations of
$\adj(X)$ when $n\leq 3$ and $K$ is a unique factorization domain. 

\begin{theorem}\label{thm:n_is_3}  
Let $K$ be a UFD and $X = (x_{ij})$ the generic $(3\times
3)$--matrix over $K$.  Then there are no nontrivial factorizations
$\adj(X)=YZ$. 
\end{theorem}

\begin{proof}  
When $n=3$, the adjoint $\adj(X)$ has determinant $(\det X)^2$, so
that a nontrivial factorization $\adj(X) = YZ$ must have $\det Y =
\det Z = \det X$.  In particular both $\cok Y$
and $\cok Z$ are of rank one and nonfree, so are each isomorphic to
one of $\{L,L^\vee\}$.  In the divisor class group $\Cl(R)$, we have $[M]
= -[L] = [L^\vee]$ by the defining sequences (\ref{eq:modulesR}), and
furthermore $[M] = [\cok Y] + [\cok Z]$ from (\ref{eq:factn-extn}).
If either of $\cok Y$ or $\cok Z$ were isomorphic to $L^\vee$, this
would force the other to be zero in $\Cl(R)$, a contradiction.  If on
the other hand both $\cok Y$ and $\cok Z$ were isomorphic to $L$, then
$[M] = 2[L]$ in $\Cl(R)$, again a contradiction. 
\end{proof}

\begin{sit}\label{sit:extensions}
Returning to the case of arbitrary $n\geq 3$ and $K$ a field, let
$\adj(X)=YZ$ be a factorization with $\det Y = \det X$.  We have the
exact sequence   
\begin{diagram}
0 & \rTo & \cok Z & \rTo & M & \rTo & \cok Y & \rTo & 0
\end{diagram} 
of Proposition~\ref{prop:corresp}, in which $\cok Y$ is isomorphic to
either $L$ or $L^\vee$.  We have also the exact sequence
\begin{diagram}
0 & \rTo & M & \rTo & \oF & \rTo & L & \rTo & 0
\end{diagram}
displaying $M$ as a first syzygy of $L$.  Form the pushout diagram:
\begin{equation}\label{eq:pushout-factn}
\begin{diagram}[midshaft]
& & & & 0 & & 0 \\
& & & & \dTo & & \dTo \\
0 & \rTo & \cok Z & \rTo & M  & \rTo & \cok Y & \rTo & 0 \\
& & \dEqual & & \dTo & \NWpbk \quad & \dTo \\
0 & \rTo & \cok Z & \rTo & \oF & \rTo & Q  & \rTo & 0 \\
& & & & \dTo & & \dTo \\
& & & & L & \rEqual & L \\
& & & & \dTo & & \dTo\\
& & & & 0 & & 0
\end{diagram}
\end{equation}
The rightmost column is an exact sequence
\begin{equation}\label{eq:our-extn}
\begin{diagram}
0 & \rTo & \cok Y & \rTo & Q & \rTo & L & \rTo & 0\,,
\end{diagram}
\end{equation}
thus naturally gives an element of $\Ext_R^1(L, \cok Y)$.
\end{sit}

\begin{prop}\label{prop:our-extn}
Let $\adj(X)=YZ$ be a nontrivial factorization of $\adj(X)$ with $\det
Y = \det X$.  Then the exact sequence (\ref{eq:our-extn}) is nonsplit,
and the middle term $Q$ is a MCM $R$-module of rank $2$ requiring at
most $n$ generators.
\end{prop}

\begin{proof}
If (\ref{eq:our-extn}) splits, then $Q\cong L \oplus \cok Y$.
Localize at the maximal ideal $\m = (x_{ij})$.  Then $\syz_1(L_\m)
\cong M_\m$ is isomorphic to a direct summand of $(\cok Z)_\m \oplus
H$ for some free $R_\m$-module $H$.  Since $M_\m$ is indecomposable
and nonfree, this implies (after passing to the completion to use the
Krull-Schmidt theorem) that $M_\m$ is a direct summand of $(\cok
Z)_\m$, and in particular that $\cok Z$ has rank at least $n-1$,
contradicting the nontriviality of the factorization $\adj(X) = YZ$. 
The statements about $Q$ follow from the diagram
(\ref{eq:pushout-factn}). 
\end{proof}

\begin{sit}
In order to classify factorizations of $\adj(X)$,
Proposition~\ref{prop:our-extn} hints that we should classify certain
extensions in $\Ext_R^1(L,\cok Y)$.  Since we assume $\det Y = \det
X$, we have either $\cok Y \cong L$ or $\cok Y \cong L^\vee$, so we must
consider $\Ext_R^1(L,L)$ and $\Ext_R^1(L,L^\vee)$.  Specifically, we are
concerned with extensions whose middle terms need the minimum number
of generators, $n$.  We define such extensions in more generality.
\end{sit}

\begin{defn}
Let $A$ be a (commutative, Noetherian) ring and $N_1$, $N_2$ finitely
generated $A$-modules.  Let $\Ext_A^1(N_1,N_2)_{\min}$ be the subset
of $\Ext_A^1(N_1,N_2)$ consisting of equivalence classes of extensions
\begin{diagram}
0 & \rTo & N_2 & \rTo & E & \rTo & N_1 & \rTo & 0
\end{diagram}
in which $E$ requires no more generators than $N_1$.
\end{defn}

\begin{lemma}\label{lem:minexts-surjs}
For any factorization $\adj(X)=YZ$ with $\det Y=\det X$, 
the natural epimorphism $\Hom_R(M,\cok Y) \to \Ext_R^1(L,\cok Y)$
induces a $1-1$ correspondence between the elements of 
$\Ext_R^1(L,\cok Y)_{\min}$ and surjective homomorphisms $M \to \cok
Y$. 
\end{lemma}

\begin{proof}
We have already seen in Proposition~\ref{prop:our-extn} that an
epimorphism $\phi \in \Hom_R(M,\cok Y)$ gives rise to an element of
$\Ext_R^1(L,\cok Y)_{\min}$.  Conversely, given an extension 
\begin{diagram}
0 & \rTo & \cok Y & \rTo & Q & \rTo & L & \rTo & 0
\end{diagram}
of $\cok Y$ by $L$, with $Q$ generated by at most $n$ elements, we
construct a commutative diagram of MCM $R$-modules with exact rows
\begin{diagram}
0 & \rTo & M' & \rTo & R^n & \rTo & L & \rTo & 0\\
&& \dTo && \dTo && \dEqual\\
0 & \rTo & \cok Y & \rTo & Q & \rTo & L & \rTo & 0\,.
\end{diagram}
By Theorem~\ref{thm:eisenbud}, there exists a pair of $(n \times
n)$--matrices $(U,V)$ giving a matrix factorization of $\det X$, and
so that $\cok U \cong L$, $\cok V \cong M'$.  It follows that $U$ is
matrix-equivalent to the generic matrix $X$, whence $V$ is equivalent
to $\adj(X)$.  Thus $M'\cong M$.  The homomorphism $M' \to \cok Y$ is
then surjective by the Snake Lemma. 

If a surjection $\phi \in \Hom_R(M,\cok Y)$ maps to zero in
$\Ext_R^1(L,\cok Y)$, that is, gives a split-exact sequence, then $Q
\cong L \oplus \cok Y$ requires more than $n$ generators.  Thus the
map from $\Hom_R(M,\cok Y)$ to $\Ext_R^1(L,\cok Y)$ is injective on
surjective homomorphisms.
\end{proof}

The results of this section set out a
correspondence between factorizations $\adj(X)=YZ$  with $\det Y =
\det X$, surjective homomorphisms $M \to \cok Y$, and elements of
$\Ext_R^1(L,\cok Y)_{\min}$.  Since $\det Y = \det X$ implies either
$\cok Y \cong L$ or $\cok Y \cong L^\vee$, we treat the two cases
separately.  In Section~\ref{section:ext1-L-L} we shall show that in fact
$\Ext_R^1(L,L)=0$, so there are {\em no} factorizations of the adjoint
with $\cok Y \cong L$.  The sections thereafter treat the case $\cok Y
\cong L^\vee$.

\section{The Module $\Ext_R^1(L,L)$ vanishes}\label{section:ext1-L-L} 

In this section we show that $L = \cok X$ has
no self-extensions, equivalently, is \emph{rigid\/} over $R=S/(\det
X)$.  Our goal follows from a recent result of R.~Ile  
\cite{Ile:2004}.  We include a proof of Ile's theorem here, since it
is short and elegant, and the matrix equations of
Section~\ref{section:matrix} simplify the argument slightly.  In the
interest of broader applicability, we will state the main results in
terms of general matrix factorizations $(\phi,\psi)$ over Noetherian
rings $S$, as Ile does, indicating where ``specialization to the
generic case'' simplifies the arguments still further.

Ile's result is couched in terms of the {\em Scandinavian complex\/}
$\Sc(\phi)$ 
attached to a matrix factorization $(\phi,\psi)$ by T.~Gulliksen and
O.~Neg\aa rd \cite{Gulliksen-Negard:1972}, which we shall have reason
to use again in Section~\ref{section:L-Lstar}.

\begin{defn}\label{def:scand}
Let $\phi: G \to F$ be a homomorphism of free modules of the same
(finite) rank $n$ over a Noetherian ring $S$.  Assume that $f = \det
\phi$ is an irreducible nonzerodivisor of $S$ and that $\ann_S\cok\phi =
(f)$, so that $(\phi,\adj(\phi),F,G)$ is a matrix factorization of
$f$.  The \emph{Scandinavian complex\/} $\Sc(\phi)$ is
\begin{equation*}\label{eq:scandcpx}
\begin{diagram}[midshaft]
0 &\rTo & S & \rTo^{?\adj(\phi)} & \Hom_S(F,G) & \rTo^{(\phi?,?\phi)} 
& \bbbh
& \rTo^{?\phi-\phi?} & \Hom_S(G,F) & \rTo^{\tr(?\adj(\phi))} & S & \rTo&
0\,,  
\end{diagram}
\end{equation*}
where $\bbbh$ is the homology in the middle of the short complex
\begin{equation*}
\begin{diagram}[midshaft] 
S & \rTo^{\Delta} & \End_S(F) \oplus \End_S(G) &  \rTo^{\tr(?)-\tr(?)}
& 
S \,,
\end{diagram}
\end{equation*}
$\tr(?)$ denotes the trace function, and $\Delta$ is the diagonal map.
\end{defn}

The complex $\Sc(\phi)$ is functorial with respect to homomorphisms of
matrix factorizations. Here is the main theorem of
\cite{Gulliksen-Negard:1972}. 

\begin{prop}[{\cite{Gulliksen-Negard:1972}; see also
	  \cite{Bruns-Vetter:1988}}]\label{prop:GNmain}
For $\phi$ as above, we have 
\[
H_0(\Sc(\phi)) \cong S/I_1(\adj(\phi)) = S/I_{n-1}(\phi)\,,
\]
and
\[
\max \{\; q \;\mid\; H_q({\Sc(\phi)}) \neq 0 \} = 4 - \grade
I_{n-1}(\phi)\,.
\]
In particular, if the grade of $I_{n-1}(\phi)$ on $S$ is $4$, the
maximum possible value, then $\Sc(\phi)$ is a (minimal, in case $S$ is
local or graded and no entry of $\phi$ is a unit) $S$-free resolution
of $S/I_{n-1}(\phi)$.
\end{prop}

\begin{remark}
It's well-known (see, for example, \cite{Hochster-Eagon:1971} or
\cite[7.3.1]{BH}) that for $\phi=X$ a generic square matrix of
indeterminates, the maximum value, $\grade I_{n-1}(X)=4$, is
achieved. 
\end{remark}

Ile's main result identifies the deformation theory of $\cok \phi$ as
the homology of $\Sc(\phi)$.  The rest of this section will be devoted
to the proof of this theorem.

\begin{theorem}[\cite{Ile:2004}]\label{thm:Ilemain}
Let $S$ be a Noetherian ring and $\phi: G \to F$ a homomorphism
between free $S$-modules of the same rank $n$, such that $f=\det \phi$
is an irreducible nonzerodivisor of $S$.  Set $R:=S/(f)$ and $M:=\cok
\phi$, and assume that $\ann_S M= (f)$.  Then 
\[
H_1(\Sc(\phi)) \cong \Ext_R^1(M,M)\,.
\]
In particular, if $\grade I_{n-1}(\phi)=4$, then $\Ext_R^1(M,M)=0$;
thus $M$ is rigid.
\end{theorem}

\begin{proof}
To compute the homology of $\Sc(\phi)$ at $\Hom_S(G,F)$, consider the
diagram
\begin{equation}\label{eq:Ilediagram}
\begin{diagram}[midshaft]
& & \bbbh & \rTo^{?\phi-\phi?} & \Hom_S(G,F) & \rTo^{\tr(?\adj(\phi))} &
  S\\ 
&& && \dTo<\pi && \dTo>{\overline{\cdot \id_M}} \\
& & & & \Ext_S^1(M,M) & \rTo_\epsilon & \End_R(M)\,.
\end{diagram}
\end{equation}
Here $\pi$ is the natural surjection and $\epsilon$ is defined by
pulling back cocycles along $\adj(\phi)$.  That is, for $\chi \in
\Ext_S^1(M,M)$, we choose a preimage $U \in \Hom_S(G,F)$ and observe that
$(U\adj(\phi),\adj(\phi)U)$ is a homomorphism of matrix factorizations
\[
(U\adj(\phi),\adj(\phi)\;U): (\phi,\adj(\phi)) \to (\phi,\adj(\phi))\,;
\]
put $\epsilon(\chi)= \cok (U\adj(\phi),\adj(\phi)\;U) \in
\End_R(M)$.  

We claim first that the square commutes.  The following lemma is the
crux of the argument.

\begin{lemma}\label{lem:Ilecrux}
For each $U \in \Hom_S(G,F)$, there exists $V \in \Hom_S(F,G)$ such
that 
\[
U\adj(\phi) - \phi V = \tr(U\adj(\phi))\cdot \id_F\,.
\]
\end{lemma}

\begin{proof} 
For the purposes of this proof, we may revert to the generic
situation, where $\phi = X$ is a square matrix of indeterminates.  As
in Section~\ref{section:matrix}, let $\partial_{ij} =
\frac{\partial}{\partial x_{ij}}$ be the partial derivative with
respect to the variable $x_{ij}$; then
\begin{align*}
\partial_{ij}[(\det X)\cdot \id_F] &= \partial_{ij}[X\adj(X)]\\
&= \partial_{ij}(X) \adj(X) + X \partial_{ij}(\adj(X))\,,
\end{align*}
where we apply $\partial_{ij}$ to a matrix entry-by-entry.  By
Lemma~\ref{lem:derdet}(1), this can be rewritten as 
\begin{equation}\label{eq:Ilecrux}
E_{ij}\adj(X) + X \partial_{ij}(\adj(X)) = \adj(X)_{ji}\cdot
\id_F\,.
\end{equation}
Write $U = (u_{ij})$; multiplying (\ref{eq:Ilecrux}) by
$u_{ij}$ and taking the sum over all $(i,j)$ gives
$$U \adj(X) - XV = \left(\sum_{i,j} u_{ij} \adj(X)_{ji}\right)\cdot
\id_F\,,$$
with $V = -\sum_{i,j} u_{ij} \partial_{ij}(\adj(X))$.  The right-hand
side of this last equation is equal to $\tr(U\adj(X))\cdot \id_F$.
\end{proof}

Returning to the proof of Theorem~\ref{thm:Ilemain}, we must show that
\[
\overline{\tr(U\adj(\phi))\cdot \id_M} = \cok(U\adj(\phi),\adj(\phi)\;U)\,,
\]
as endomorphisms of $M$, for each $U \in \Hom_S(G,F)$.  By the Lemma,
there exists $V \in \Hom_S(F,G)$ so that 
\[
U\adj(\phi) - \phi V = \tr(U\adj(\phi))\cdot \id_F\,.
\]
In particular, the two sides induce the same endomorphism of $M$.  The
term $\phi V$ factors through $F$, so gives the zero map on $M = \cok
\phi$; thus $U\adj(\phi)$ induces $\overline{\tr(U\adj(\phi))\cdot
  \id_M}$.

Next we shall show that $\ker \epsilon \cong \Ext_R^1(M,M)$.  Indeed,
an $S$-module extension $\chi$, represented by $U \in \Hom_S(G,F)$,
is an extension of $R$-modules if and only if $U$ is part of a
homomorphism of matrix factorizations, {\it i.e.\/}, there exists $V
\in \Hom_S(F,G)$ so that $U\adj(\phi) = \phi V$.  This is the case
precisely when $U\adj(\phi)$ factors through $G$, that is, induces the
zero endomorphism of $M$. 

Finally, we claim that $\pi$ induces an isomorphism $H_1(\Sc(\phi))
\to \ker \epsilon$.  To see this, first let $[U]$ be a homology
class.  Then the image of $U$ in $\End_R(M)$ is zero, so that $\pi(U)
\in \ker \epsilon$.  Next, take $U \in \Hom_S(G,F)$ to be a boundary,
so that $U = A\phi-\phi B$ for some $(A,B) \in \bbbh$.  Then the
homomorphism of matrix factorizations induced by $\pi(U)$ is
equivalent to $(\phi B \adj(\phi), \adj(\phi) A \phi)$.  Since $\phi B
\adj(\phi)$ factors through $G$, this is zero in $\End_R(M)$.  Lastly,
any $\chi \in \Ext_R^1(M,M)$ lifts to $U \in \Hom_S(G,F)$, which must
then be a cycle by the commutativity of the square.  This finishes the proof.
\end{proof}

Specializing to the case of a generic matrix, we obtain our main
result of this section. 

\begin{cor}\label{cor:mainL-L}
Let $K$ be a commutative Noetherian normal domain, $X=(x_{ij})$ the
generic $(n \times n)$--matrix over $K$, $S = K[x_{ij}]$, and $R =
S/(\det X)$.  Set $L:=\cok X$.  Then 
\begin{enumerate}
\item $\End_R(L) \cong R$\,;
\item $\Ext_R^1(L,L) = 0$\,;
\item $\Ext_S^1(L,L)$ is isomorphic to the ideal $I_{n-1}(X)/(\det X)$
  of $R$\,;
\item $\Ext_R^2(L,L) \cong S/I_{n-1}(X)$\,; and
\item $\uEnd_R(L) \cong S/I_{n-1}(X)$\,.
\end{enumerate}
\end{cor}

\begin{proof}
We have already observed that $R$ is a normal domain (see
\ref{notation}).  Since 
$L$ has rank one, the ring $\End_R(L)$ is a finite extension of $R$
contained in its quotient field, so equal to $R$ by normality.
Claims (2), (3), and (4) follow from
Theorem~\ref{thm:Ilemain} and the diagram (\ref{eq:Ilediagram}): Since
$\grade I_{n-1}(X)=4$, we have $\Ext_R^1(L,L)=0$, and  the
image of $\epsilon$ is equal to the image of $\tr(?\adj(X))$, that is,
$I_{n-1}(X)$.  Finally, statement (5) follows from (4) and
Proposition~\ref{prop:stable-ext-periodic}. 
\end{proof}

\section{The module $\Hom_R(M,L^\vee)$ and Classification of
  Factorizations} \label{section:hom-M-Lstar}

In this section we consider the case $\cok Y \cong
L^\vee$.  Using the matrix equations of Section~\ref{section:matrix},
we compute the 
free resolution of the $S$-module $\Hom_R(M,L^\vee)$, and compare the
result to a canonical short exact sequence 
 to show that not only does every
factorization $\adj(X)=YZ$ with $\det Y = \det X$ yield an extension
of $L$ by $L^\vee$, but we can classify precisely when two factorizations
give equivalent extensions.

Our next task is to interpret the matrix-theoretic result
Theorem~\ref{thm:matrixmain} in terms of MCM modules over $R$.  Keep
the notation of \ref{notation}.

\begin{remark}\label{rmk:hom-def}
Let $A \in \Alt_n(S)$ be an alternating $(n\times n)$--matrix over $S$
and let $B_A$ be the companion matrix of Theorem~\ref{thm:matrixmain},
so that 
\[
A \adj(X) = X^T B_A\,.
\]
This equation defines a commutative diagram of free $S$-modules
\begin{equation*}
\begin{diagram}[midshaft]
G(-n) & \rTo^{X} & F(-n) & \rTo^{\adj(X)} & G \\
\dTo<{A} & & \dTo<{B_A} & & \dTo>{A} \\
F(-n) & \rTo_{\adj(X)^T} & G & \rTo_{X^T} & F\,,
\end{diagram}
\end{equation*}
that is, a homomorphism of matrix factorizations
\[
(A,B_A): (\adj(X),X) \to (X^T,\adj(X)^T)
\]
and thus a homomorphism of MCM $R$-modules 
\[
\cok(A,B_A) : M \to L^\vee.
\]
In other words, we have a homomorphism $\Alt_n(S) \xto{A \mapsto
  \cok(A,B_A)} \Hom_R(M,L^\vee)$.  Our next result is that this
  homomorphism is surjective, so that $\Hom_R(M,L^\vee)$ is generated
  by the alternating matrices, and moreover that $\Hom_R(M,L^\vee)$ is
  itself a MCM $R$-module.
\end{remark}

\begin{theorem}\label{thm:hom-res}
The $R$-module $\Hom_R(M,L^\vee)$ is MCM of rank $n-1$, minimally
generated by $\binom n 2$ elements.  More precisely, it has the
following free presentation as an $S$-module
\begin{diagram}[midshaft]
0 & \rTo& \Alt_n(S) & \rTo^{U \mapsto X^T U X} & \Alt_n(S)(2) & \rTo^{A
  \mapsto \cok(A,B_A)} & \Hom_R(M,L^\vee) & \rTo & 0 \,;
\end{diagram}
alternatively, in terms of exterior powers, this exact sequence
can be written as  
\begin{diagram}[midshaft]
0 & \rTo & \Lambda^2 F^\vee & \rTo^{\Lambda^2 X^T}& \Lambda^2 G^\vee & \rTo
& \Hom_R(M,L^\vee) & \rTo & 0\,.
\end{diagram}
\end{theorem}

\begin{proof}
For $U$ an alternating $(n \times n)$--matrix over $S$, we have
\begin{align*}
X^TB_{X^TUX} &= (X^TUX)\adj(X) \\
   &= X^T U \cdot (\det X)\,,
\end{align*}
so that $B_{X^TUX} = U\cdot (\det X)$.  Thus the homomorphism
$\cok(X^TUX,B_{X^TUX})$ is zero on the $R$-module $M$, and the alleged
resolution of $\Hom_R(M,L^\vee)$ is at least a complex.

Put $D := \cok (U \mapsto X^TUX)$.  Then $D$ maps to $\Hom_R(M,L^\vee)$
and we must show that this map is an isomorphism.  Note first that $D$
is a MCM $R$-module, with matrix factorization 
\[
(U \mapsto X^TUX, A \mapsto B_A)\,.
\]
Indeed, we have seen that $B_{X^TUX} = U \cdot (\det X)$, and also
$X^TB_AX = A \adj(X) X = A \cdot (\det X)$.  Thus in particular $D$ is
a reflexive $R$-module, and $U \mapsto X^TUX$ is an injective
endomorphism of the module of alternating matrices.

The free module $\Alt_n(S)$ has rank $\binom n 2$, so the determinant
of the endomorphism $U \mapsto X^TUX$ is homogeneous of degree
$n(n-1)$ in the variables $x_{ij}$.  Since it must also be a unit
times $(\det X)^{\rank D}$, we see that $D$ has rank $n-1$ as an
$R$-module, equal to that of $\Hom_R(M,L^\vee)$.

Outside the singular locus $V(I_{n-1}(X))$ of $R$, at least one
maximal minor of $X^T$ is a unit.  Thus after elementary transformations
and linear changes of variables, $X^T = \diag(0, 1, \dots, 1)$ and so
$\adj(X)^T = E_{11}$, the elementary matrix with $1$ at position $(1,1)$
and zeros elsewhere.  Now any homomorphism $\alpha$ from the cokernel
of $E_{11}$ to the cokernel of $\diag(0,1,\dots,1)$ is induced by an
alternating $(n \times n)$--matrix, namely any alternating matrix with
first row $\alpha$.  That is, outside the singular locus of $R$,
$\Hom_R(M,L^\vee)$ is indeed generated by homomorphisms $\cok(A,B_A)$ for
alternating $A$.  The map $D \to \Hom_R(M,L^\vee)$ is thus surjective,
and since $D$ and $\Hom_R(M,L^\vee)$ have the same rank, is even an
isomorphism, outside $V(I_{n-1}(X))$.

Recall that $R$ is normal and $I_{n-1}(X)$ has codimension $3$ in
$\Spec R$.  The homomorphism $D \to \Hom_R(M,L^\vee)$
is thus a homomorphism between reflexive modules over a normal domain,
which is an isomorphism in codimension one.  It follows that in fact
$D \to \Hom_R(M,L^\vee)$ is an isomorphism.
\end{proof}

\begin{sit}\label{sit:other-hom-res}
Define an $S$-module $C$ by the pullback diagram
\begin{diagram}[midshaft]
C & \rTo && \Hom_S(G,G^\vee) \\
\dTo & & & \dTo>{?\adj(X)} \\
\Hom_S(F(-n),F^\vee) & \rTo_{X^T?} && \Hom_S(F(-n),G^\vee)
\end{diagram}
Then 
\[
C = \{ (A,B) \;|\; A \adj(X) = X^T B\}\,.
\]
By Remark~\ref{rmk:hom-def}, there is a natural homomorphism
$C \to \Hom_R(M,L^\vee)$, sending $(A, B)$ to $\cok(A,B)$.  There is also
a natural embedding 
$\Hom_S(G,F^\vee) \to C$, given by $U \mapsto (X^TU,U\adj(X))$.  This
gives an exact sequence of $S$-modules ({\em cf.\/}
\ref{sit:matfacthoms})
\begin{equation}\label{eq:other-res-hom}
\begin{diagram}
0 & \rTo & \Hom_S(G,F^\vee) & \rTo & C & \rTo & \Hom_R(M,L^\vee) & \rTo & 0\,.
\end{diagram}
\end{equation}
Comparing the resolution given by Theorem~\ref{thm:hom-res} to
(\ref{eq:other-res-hom}) will give our classification of factorizations
of $\adj(X)$.  To prepare for this, we make a definition.
\end{sit}

\begin{defn}
A factorization of the form $\adj(X) = JXZ$ or $\adj(X)=JX^TZ$, with
$J$ an invertible $(n \times n)$--matrix, is called a {\em
  normalized factorization\/} of the adjoint.
\end{defn}

Note that since, up to equivalence, the only $(n\times n)$--matrices
with determinant equal to $\det X$ are $X$ and $X^T$, we may always
normalize a given factorization.  Explicitly, if $\adj(X)=YZ$ with
$\cok Y \cong L^\vee$, there exist invertible matrices $J$ and $J'$ so
that $Y=JX^TJ'$, and replacing $Z$ by $J'Z$ gives a normalized
factorization.  Similarly, if $\cok Y \cong L$, we may multiply $Z$ on
the left by an invertible matrix to achieve a normalized factorization.

\begin{theorem}\label{thm:normfactclassify}
Let $\adj(X) = YZ = JX^TZ$ be a normalized factorization of $\adj(X)$
with $\cok Y \cong L^\vee$.  Then there exist a unique invertible
alternating $(n \times n)$--matrix $A$ and a unique $(n \times
n)$--matrix $U$ such that  
\[
J^{-1} = A + X^TU {\text{ \quad and \quad }} Z = B_A + U\adj(X)\,.
\]
In particular, the homomorphism $\cok(A,B_A) : M \to L^\vee$ is
surjective.  Two normalized factorizations $JX^TZ$ and $J'X^TZ'$ give
the same epimorphism in $\Hom_R(M,L^\vee)$ if and only if $J^{-1} -
J'^{-1} = X^TV$ for some $(n \times
n)$--matrix $V$, and then $Z-Z' = V\adj(X)$. 
\end{theorem}

\begin{proof}
We have a commutative diagram of $S$-modules with exact rows
\begin{diagram}[midshaft]
0 & \rTo & \Lambda^2 F^\vee & \rTo^{X^T?X} & \Lambda^2 G^\vee & \rTo &
\Hom_R(M,L^\vee) & \rTo & 0\\
& & \dTo<{?X} & & \dTo>{A \mapsto (A,B_A)} && \dEqual\\
0 & \rTo & \Hom_S(G,F^\vee) & \rTo_{(X^T?,?\adj(X))}& C & \rTo &
\Hom_R(M,L^\vee) & \rTo & 0
\end{diagram}
in which the vertical arrows represent monomorphisms.
The factorization $\adj(X)=JX^TZ$ yields
$J^{-1}\adj(X) = X^TZ$,
so gives a homomorphism $(J^{-1},Z)$ of matrix factorizations.  From
the diagram, we obtain $A \in \Lambda^2 G^\vee$, unique by the Snake
Lemma, so that $J^{-1} = A + X^T U$ for some $U \in \Hom_S(G,F^\vee)$ and
$Z = B_A + U \adj(X)$.  Since $(J^{-1},Z)$ and $(A,B_A)$ induce the
same homomorphism $M \to L^\vee$, $A$ is also invertible.  The final
statement follows from the uniqueness of $A$ and \ref{sit:matfacthoms}
above.
\end{proof}

\begin{cor}\label{cor:correspondence}
There is a one-to-one correspondence between equivalence classes of
extensions 
\begin{diagram}
0 & \rTo & L^\vee & \rTo & Q & \rTo & L & \rTo & 0
\end{diagram} 
in which $Q$ is a homomorphic image of $R^n$ and equivalence classes
of normalized factorizations $\adj(X) = JX^TZ$, where the equivalence
relation is given by 
\[
JX^TZ \sim J'X^TZ' \quad \text{ iff } \quad J^{-1}-J'^{-1} = X^TV
\text{ and } Z-Z'= V\adj(X)
\]
for some $(n \times n)$--matrix $V$.
\end{cor}

\section{The Module $\Ext_R^1(L,L^\vee)$}\label{section:L-Lstar}

We now turn to computing $E := \Ext_R^1(L,L^\vee)$.  
Keep the notation of \ref{notation}.

The epimorphism $\Lambda^2 G^\vee \to \Hom_R(M,L^\vee)$ of
Theorem~\ref{thm:hom-res} composes with the natural epimorphism
$\Hom_R(M,L^\vee) \to E$ to begin a free resolution of $E$.  Let us first
identify this map more explicitly.  Recall that a finitely generated
module $N$ over a normal domain is {\em orientable\/} provided $[N]=0$
in the divisor class group.

\begin{lemma}\label{lem:rk2exist}
For each alternating $(n\times n)$--matrix $A$ over $S$, there exists
an extension
\begin{equation}\label{eq:Qext}
\begin{diagram}
0 & \rTo & L^\vee & \rTo & Q & \rTo & L & \rTo &0\,,
\end{diagram}
\end{equation}
which is the image of $\cok(A,B_A)$ under the natural epimorphism
$\Hom_R(M,L^\vee) \to \Ext_R^1(L,L^\vee)$.  In particular, the module $Q$ is
a MCM $R$-module of rank $2$ given by the matrix factorization
\[
Q = \cok\left(
\left(\begin{matrix} X^T & A \\ 0 & X \end{matrix}\right), 
\left(\begin{matrix} \adj(X)^T & -B_A \\ 0 & \adj(X)
\end{matrix}\right) 
\right)\,.
\]
Furthermore, $Q$ is orientable.
\end{lemma}

\begin{proof}
This is a restatement of Remark~\ref{rmk:hom-def} and
Proposition~\ref{sit:pushout}. To see that $Q$ is orientable,
observe that $[Q] = [L]+[L^\vee]=[L]-[L] =0$ in $\Cl(R)$.
\end{proof}

\begin{remark}
The matrix factorization given for $Q$ in Lemma~\ref{lem:rk2exist} may
not be of minimal 
size.  Indeed, if $A$ is invertible then we have seen that $Q$
requires only $n$ generators.  In this case, we have
\[
\left(\begin{matrix} \id_n & 0 \\ -XA^{-1} & \id_n \end{matrix}\right)
\left(\begin{matrix} X^T & A \\ 0 & X \end{matrix}\right)
\left(\begin{matrix} \id_n & 0 \\ -A^{-1}X^T & \id_n \end{matrix}\right)
=
\left(\begin{matrix} 0 & A \\ -XA^{-1}X^T & 0 \end{matrix}\right)\,,
\]
so that the given matrix
factorization for $Q$ can be reduced to 
\[
Q \cong \cok(XA^{-1}X^T,B_A)\,.
\]
More generally, if a $(k\times k)$--minor of $A$ is invertible, then
the given matrix factorization of $Q$ can be reduced to one of size
$2n-k$.

For any $n$, the graded, orientable rank two MCM $R$-modules are
minimally evenly generated \cite[3.1]{Herzog-Kuhl:1989}.  In fact,
they are presented by yet another alternating matrix over $S$, as in
\cite{Baciu-Ene-Pfister-Popescu:2004}. 
\end{remark}

\begin{remark}
The orientable MCM module $Q$ of Lemma~\ref{lem:rk2exist} is
decomposable if and only if $Q \cong L \oplus L^\vee$, equivalently,
the sequence (\ref{eq:Qext}) is split exact.  To see this, recall that
$L$ and $L^\vee$ are up to isomorphism the only MCM $R$-modules of rank
one. As $Q$ is orientable, the only possible direct-sum decomposition
for $Q$ is $L \oplus L^\vee$, and by Miyata's theorem \cite{Miyata}, if
(\ref{eq:Qext}) is apparently split then it is split.
\end{remark}

We first determine the structure of $E$ as an $S$-module.

\begin{theorem}\label{thm:res-E}
The $S$-module $E = \Ext_R^1(L,L^\vee)$ has the following graded minimal
resolution.
\begin{equation*}\label{eq:Extres}
\begin{diagram}[small,midshaft] 
&&&&&\bbbs_{2}G^\vee(-n)\\
&&&&\ruTo^{X^T?+?^TX}&&\rdTo^{-\adj(X)}\\
0\rightarrow &\Lambda^{2}F^\vee(-n) & \rTo^{\ ?X} & F^\vee\otimes
G^\vee(-n) &&&& G^\vee\otimes F^\vee&
\rTo^{?X-X^T?^T}&\Lambda^{2}G^\vee &\rTo & E \rightarrow 0\\
&&&&\rdTo_{?\adj(X) + \adj(X)^T?^T}&&\ruTo_{X^T?}\\
&&&&&\bbbd_{2}F^\vee&&
\end{diagram}
\end{equation*}
\end{theorem}

\begin{proof}
The diagram of $S$-modules
\small
\begin{displaymath}
\begin{diagram}[midshaft,width=2em]
& && && && 0 && \\
& && && && \uTo \\
& && && && \Ext_R^1(L,L^\vee) \\
& && && && \uTo \\
& 0 & \rTo & \Lambda^2 F^\vee & \rTo^{X^T?X} & \Lambda^2 G^\vee & \rTo &
  \Hom_R(M,L^\vee) 
  & \to  0  \\
& && \uTo<{?-?^T} && \uTo>{?X-X^T?^T} && \uTo  \\
& 0 & \rTo & \Hom_S(F,F^\vee) & \rTo^{X^T?} & \Hom_S(F,G^\vee) & \rTo &
  \Hom_R(\oF,L^\vee) & \to  0 \\
& && \uTo<{?\adj(X) + \adj(X)^T?^T} && \uTo>{?\adj(X)} && \uTo \\
0 \to & \Lambda^2F^\vee(-n) & \rTo_{\ ?X} & F^\vee \otimes_S G^\vee(-n) &
\rTo_{X^T?+?^TX\ \ } & \bbbs_2G^\vee(-n) & \rTo & \Hom_R(L,L^\vee) & \to  0 \\
& && && && \uTo \\
& && && && 0
\end{diagram}
\end{displaymath}
\normalsize
commutes, has exact rows, and has complexes for columns.  The
ingredients of the diagram are as follows.  For rows, it has the
resolution of $\Hom_R(M,L^\vee)$ computed in
Theorem~\ref{thm:hom-res}, the result of applying $\Hom_S(F,-)$
to the $S$-module resolution of $L^\vee$, and the resolution of
$\Hom_R(L,L^\vee) \cong \bbbs_2 L^\vee$ computed via the Eagon--Northcott
complex in \cite[Thm. 2.16]{Bruns-Vetter:1988}.  The rightmost column
is the result of applying $\Hom_R(-,L^\vee)$ to the $R$-module resolution
of $L$.  The other columns are liftings of the maps in the rightmost
column. 

Truncate the diagram, retaining only the free $S$-modules; since the
rows are acyclic, the total complex is acyclic as well, with zeroth
homology isomorphic to $E$.  This gives a free resolution
\small
\begin{displaymath}
\begin{diagram}[height=.85em,midshaft] 
&&&& G^\vee\otimes F^\vee && \bbbs_{2}G^\vee(-n) &\\
0\leftarrow E &\lTo&\Lambda^{2}G^\vee&\lTo& {\qquad\oplus\qquad}  &\lTo &
{\qquad\oplus\qquad}  &\lTo 
&F^\vee\otimes G^\vee(-n)& \lTo&\Lambda^{2}F^\vee(-n) \leftarrow  0 \\
&&&& \Lambda^2 F^\vee && F^\vee \otimes F^\vee &\\
\end{diagram}
\end{displaymath}
\normalsize
of $E$.  By \ref{sit:alt-symm-triang}, the map $F^\vee \otimes F^\vee
\xto{?-?^T} \Lambda^2 F^\vee$ splits 
out a direct summand isomorphic to the kernel, $\bbbd_2F^\vee$, which
gives the resolution claimed. 
\end{proof}

\begin{cor}\label{cor:Hilb-series-E}
As a graded module, the Hilbert series of $E$ is
$$
H_{E}(t) =\frac{{\binom n 2}t^2-n^{2}t^3+{\binom {n+1}{2}
}(t^{4}+t^{n+2})-n^{2}t^{n+3}+{\binom n 2}t^{n+4}}
{(1-t)^{n^{2}}}\,.
$$
\end{cor}

\begin{prop}\label{prop:supp-E}
As an $S$-module, $E$ is perfect of grade $4$, with support the
singular locus $V(I_{n-1}(X))$ of $R$.  More precisely, the
annihilator of $E$ is equal to $I_{n-1}(X)$.
\end{prop}

\begin{proof}
To see that $I_{n-1}(X)$ annihilates $E$, fix an index $i$ and denote
by $\bar X_i$ the $n\times(n-1)$--matrix obtained by deleting the
$i^\text{th}$ column of $X$.  Then $L^\vee$ is isomorphic to
$I_{n-1}(\bar X_i)$ for each $i$.  The natural epimorphism
$\Hom_R(L,R/I_{n-1}(\bar X_i)) \to E$ shows that $I_{n-1}(\bar X_i) E
=0$.  On the other hand, $L$ is isomorphic to $I_{n-1}(\bar X_i')$,
where $\bar X_i'$ is obtained by deleting the $i^\text{th}$ {\em
  row\/} of $X$.  Thus $E \cong \Ext_R^2(R/I_{n-1}(\bar X_i'),L^\vee)$
is also killed by $I_{n-1}(\bar X_i')$.  Letting $i$ vary, we see that
$E$ is annihilated by every $(n-1)\times(n-1)$--minor of $X$, that is,
by $I_{n-1}(X)$. 

By the resolution of Theorem~\ref{thm:res-E}, $\ann_S E$ has
codimension at most
$4$.  But $I_{n-1}(X) \subseteq \ann_S E$ and $I_{n-1}(X)=4$ is a
prime ideal of height $4$, so that $I_{n-1}(X) = \ann_S E$.
\end{proof}

It follows from Proposition~\ref{prop:supp-E} that $E$ is naturally a
module over $S/I_{n-1}(X)$.  The next result details the structure of
$E$ as an $S/I_{n-1}(X)$-module.  

\begin{theorem}\label{thm:ext-mcm-ideal}
As an $S/I_{n-1}(X)$-module, $E = \Ext_R^1(L,L^\vee)$ is a MCM module of
rank one, isomorphic to the ideal generated by the maximal minors of
$n-2$ fixed rows of $X$.
\end{theorem}

\begin{proof}
Fix $r,s$ with $1 \leq r < s \leq n$ and define a homomorphism
$\xi_{rs}: E \to S/I_{n-1}(X)$ as follows.  For an alternating matrix
$A$, let $B_A = (b_{ij})$ be the companion matrix of
Theorem~\ref{thm:matrixmain}, so that
\[
A \adj(X) = X^T B_A\,.
\]
Set $\xi_{rs}(A) = b_{rs}$.  To verify that
$\xi_{rs}$ defines a well-defined homomorphism on $E$, it 
suffices (in view of Theorem~\ref{thm:res-E}) to show that
$\xi_{rs}(UX-X^TU^T)=0$ for any $(n \times n)$--matrix $U$.  Since the
companion $B_{UX-X^TU^T}$ satisfies
\[
(UX-X^TU^T)\adj(X) = X^TB_{UX-X^TU^T}
\]
and is the unique matrix with this property, we see that 
\[
B_{UX-X^TU^T} = \adj(X)^TU - U^T\adj(X)\,,
\]
so that $b_{rs} \in I_1(\adj(X)) = I_{n-1}(X)$ and
$\xi_{rs}(UX-X^TU^T)=0$ in $S/I_{n-1}(X)$.

Recall that 
\[
b_{rs} = \sum_{k<l} a_{kl}(-1)^{r+s+k+l}[rs\;\widehat\mid\;kl]\,.
\]
The image of $\xi_{rs}$, as $A$ ranges over all alternating matrices,
is thus equal to the ideal generated by all minors
$[rs\;\widehat\mid\;kl]$ for $1 \leq k < l \leq n$.

To show that 
$\xi_{rs}$ is a monomorphism, it suffices to show that 
$E$ has rank one over the integral domain $S/I_{n-1}(X)$, or
equivalently that we have an equality of multiplicities $e(E) =
e(S/I_{n-1}(X))$.  The (graded) Scandinavian complex
\small
\begin{diagram}
0 \to  S(-2n) & \rTo^{?\adj(X)} & G^\vee\otimes_S F(-n) & \rTo^{(X?,?X)} 
& \bbbh 
& \rTo^{?X-X?} & F^\vee\otimes_S G(-n) & \rTo^{\tr(?\adj(X))} & S \to 0 
\end{diagram}
\normalsize
of Definition~\ref{def:scand} is a graded minimal $S$-free resolution
of $S/I_{n-1}(X)$, so from it 
we compute the Hilbert series
\[
H_{S/I_{n-1}(X)}(t) = \frac{1-n^2t^{n-1}+ (2n^2-2)t^n-n^2t^{n+1} +t^{2n}}{(1-t)^{n^2}}
\]
and see from Corollary~\ref{cor:Hilb-series-E} that $e(E) =
\frac{1}{12}(n^4-n^2) = e(S/I_{n-1}(X))$ ({\em cf.\/} also
\cite{Herzog-Trung:1992}). 
\end{proof}

Consider the homomorphism of free $S$-modules
\[
B_?: \Lambda^2 G^\vee \to \Lambda^2 F^\vee
\]
that sends an alternating matrix $A$ to its companion matrix $B_A$,
satisfying $A \adj(X) = X^T B_A$.  Recall from
Theorem~\ref{thm:hom-res} that $\Hom_R(M,L^\vee)$ is isomorphic to the
image of $B_? \otimes_S R$.  Here is another characterization of
$E=\Ext_R^1(L,L^\vee)$ in these terms, which follows from the description
of the epimorphism $\Lambda^2 G^\vee \to E$ in Theorem~\ref{thm:res-E}.

\begin{prop} \label{prop:ext-image}
The module $\Ext_R^1(L,L^\vee)$ is isomorphic to the image of the
homomorphism of free $S/I_{n-1}(X)$-modules $B_? \otimes_S
S/I_{n-1}(X)$. 
\end{prop}

\begin{remark}\label{rmk:symmetry}
Of course one has symmetric results for $E' = \Ext_R^1(L^\vee,L)$,
exploiting the fact that $\tau^*L \cong L^\vee$ and so $\tau^*E \cong
E'$, where $\tau$ is again the involution of $S$ defined in
\ref{sit:eval}.  In particular, $E'$ is also a MCM module of rank one
over $R_1:= S/I_{n-1}(X)$, isomorphic to the ideal generated by the
maximal minors of any $n-2$ fixed {\em columns\/} of $X$.  It follows
that in fact $E$ and $E'$ are, up to isomorphism, the only rank-one
MCM $S/I_{n-1}(X)$-modules, opposites of each other in
$\Cl(R_1)$ \cite[9.27]{Bruns-Vetter:1988}.  
\end{remark}

\section{Extensions of rank-one MCM modules}\label{section:extsofrk1s}

The results of Section~\ref{section:L-Lstar} classify the extensions
of $L^\vee$ by $L$ up to equivalence.  Of course inequivalent extensions
may have isomorphic middle terms.  In this section we describe the MCM
modules over the generic determinantal hypersurface ring which appear
as middle terms of the extensions in Theorem~\ref{thm:res-E}.  A
complete classification seems out of reach; as soon as $n \geq 3$, the
class of such modules cannot be parametrized by the points of any
finite-dimensional algebraic variety (Corollary~\ref{cor:wild}).

First consider the case $n=2$.  We recover from
Theorem~\ref{thm:ext-mcm-ideal} the following consequence of
\cite{BGS}.

\begin{cor}\label{cor:n=2}
Let $K$ be a field and $X=(x_{ij})$ the generic $(2 \times 2)$--matrix
over $K$.  Then there are no indecomposable MCM modules of rank $2$
over $R = K[x_{ij}]/(\det X)$ which are extensions of rank-one MCM
$R$-modules. 
\end{cor}

\begin{proof}
As $\Ext_R^1(L,L^\vee) \cong K$ is a one-dimensional vector space,
there are up to equivalence only two nonsplit extensions of rank-one 
MCM $R$-modules: the short exact sequences displaying $L^\vee$ as a
first syzygy of $L$, and vice versa.  The middle terms of each of
these are free of rank two, not indecomposable.
\end{proof}

Of course, far more is true in this case: When $n=2$ and $K$ is an
algebraically closed field, the {\em only\/} indecomposable MCM
$R$-modules are $R$, $L$, and $L^\vee$ by the classification in
\cite{BGS}, so that $R$ has {\em finite Cohen--Macaulay type\/}.  This
property fails dramatically for $n\geq 3$.

To describe this failure precisely, let us suspend all our notational
assumptions for a moment and consider a more general problem.

\begin{sit}\label{sit:gen-exts-and-resns}
Let $A$ and $B$ be finitely generated modules over a (commutative,
Noetherian) ring $R$.  Fix
free resolutions
\begin{equation*}
\begin{diagram}
\cdots & \rTo & P_2 & \rTo^{X_2} & P_1 & \rTo^{X_1} & P_0 & \rTo & A & \rTo &
0 \\
\cdots & \rTo & Q_2 & \rTo^{Y_2} & Q_1 & \rTo^{Y_1} & Q_0 & \rTo & B & \rTo &
0
\end{diagram}
\end{equation*}
of $A$ and $B$.  An element $\chi \in \Ext_R^1(A,B)$ is an equivalence
class of extensions 
\begin{equation*}
\begin{diagram}
0 & \rTo & B & \rTo & E & \rTo & A & \rTo & 0\,,
\end{diagram} 
\end{equation*}
and the isomorphism class of $E$ is determined by $\chi$.  The
Horseshoe Lemma provides a free resolution of $E$
\begin{equation*}
\begin{diagram}
\cdots & \rTo & Q_2\oplus P_2 & 
\rTo^{\left[\begin{smallmatrix}Y_2 & Z_2 \\ 0 &
	  X_2\end{smallmatrix}\right]} 
& Q_1\oplus P_1 & 
\rTo^{\left[\begin{smallmatrix}Y_1 & Z_1 \\ 0 &
	  X_1\end{smallmatrix}\right]} 
& Q_0\oplus P_0 & \rTo & E & \rTo &0\,.
\end{diagram}
\end{equation*}
Here the $Z_i$ are homomorphisms in $\Hom_R(P_i,Q_{i-1})$ satisfying
$Y_iZ_{i+1}+Z_iX_{i+1} =0$ for all $i \geq 1$.
\end{sit}

\begin{defn}
In the situation of \ref{sit:gen-exts-and-resns}, define a sequence
of rings
\[
R_i := R/(I_1(X_i) + I_1(Y_i))
\]
for $i = 1, 2, \dots$, where as usual $I_1(U)$ is the ideal of $R$
generated by the entries of $U$.  For each $i$ set  
\[
\J_i(\chi) = \frac{I_1(Z_i) + I_1(X_i) + I_1(Y_i)}{I_1(X_i) +
  I_1(Y_i)}\,,
\]
an ideal of $R_i$.
\end{defn}

It is straightforward to check that the ideals $\J_i(\chi) \subseteq
R_i$ are well-defined.  In fact, the $\J_i(\chi)$ are invariants of
the isomorphism class of the middle term of $\chi$:

\begin{prop}
Let $\chi, \chi' \in \Ext_R^1(A,B)$ have middle terms $E, E'$.  If $E
\cong E'$, then $\J_i(\chi) = \J_i(\chi')$ for all $i$.
\end{prop}

The function $\J_i(?)$ thus defines a map from isomorphism classes of
modules $E$ appearing as extensions of $B$ by $A$ to ideals of
$R_i$.  We can identify which ideals are in the image of $\J_1$.

\begin{prop}
Let $Z_1: P_1 \to Q_0$ and $Z_2:P_2 \to Q_1$ be homomorphisms of free
modules such that $Y_1 Z_2 + Z_1X_2=0$.  Then there exists $\chi \in
\Ext_R^1(A,B)$ such that $\J_1(\chi) = I_1(Z_1)R_1$.
\end{prop}

\begin{proof}
Set $E = \cok \left[\begin{smallmatrix}Y_1 & Z_1 \\ 0 &
	  X_1\end{smallmatrix}\right]$, so that we have a commutative
	  diagram
\begin{equation*}
\begin{diagram}
& & 0 & & 0 & & 0 \\
& & \uTo & & \uTo & & \uTo \\
& & B & \rTo & E & \rTo & A \\
& & \uTo && \uTo && \uTo \\
0 & \rTo & Q_0 & \rTo & Q_0 \oplus P_0 & \rTo & P_0 & \rTo & 0 \\
& & \uTo<{Y_1} & & \uTo<{\left[\begin{smallmatrix}Y_1 & Z_1 \\ 0 &
	  X_1\end{smallmatrix}\right]} & & \uTo>{X_1} \\
0 & \rTo & Q_1 & \rTo & Q_1 \oplus P_1 & \rTo & P_1 & \rTo & 0 \\
& & \uTo<{Y_2} & & \uTo<{\left[\begin{smallmatrix}Y_2 & Z_2 \\ 0 &
	  X_2\end{smallmatrix}\right]} & & \uTo>{X_2} \\
0 & \rTo & Q_2 & \rTo & Q_2 \oplus P_2 & \rTo & P_2 & \rTo & 0
\end{diagram}
\end{equation*}
with exact rows and columns.  The map $E \to A$ is surjective by
commutativity.  To see that $B \to E$ is injective, it is equivalent
by the Snake Lemma to see that the kernel of
$\left[\begin{smallmatrix}Y_1 & Z_1 \\ 0 &
	X_1\end{smallmatrix}\right]$ maps onto the kernel of $X_1$.  This
is a straightforward calculation using $Y_1 Z_2 + Z_1 X_2=0$.
\end{proof} 

\begin{sit}\label{sit:gen-hyper}
Assume now that $R = S/(f)$ is a hypersurface ring and $A, B$ are MCM
modules over $R$.  The free resolutions of $A$ and $B$ are periodic of
period $2$, given by matrix factorizations of $f$.  Write $A = \cok
(\phi,\psi)$ and $B = \cok (\phi',\psi')$.  Then the sequence of rings
$R_i$ is periodic: we have
\[
R_i = 
  \begin{cases}
	S/(I_1(\phi) + I_1(\phi')) & \text{for $i$ odd, and}\\
	S/(I_1(\psi) + I_1(\psi')) & \text{for $i$ even.}
  \end{cases}
\]
For $\chi \in \Ext_R^1(A,B)$, the ideals $\J_1(\chi) \subseteq R_1$
and $\J_2(\chi) \subseteq R_2$ are again invariants of the middle term
of $\chi$.  
\end{sit}

\begin{sit} \label{ideals-gen-ext}
Return now to the generic determinant, with notation as in
\ref{notation}.  Consider $\Ext_R^1(L,L^\vee)$.  Since $L = \cok(X,
\adj(X))$ and $L^\vee = \cok(X^T, \adj(X)^T)$, we have
\[
R_i =
  \begin{cases}
	S/I_1(X) \cong K & \text{for $i$ odd, and}\\
	S/I_1(\adj(X) = S/I_{n-1}(X) & \text{for $i$ even.}
  \end{cases}
\]
By Theorem~\ref{thm:res-E} and Lemma~\ref{lem:rk2exist}, every element
$\chi \in \Ext_R^1(L,L^\vee)$ is of the form
\begin{equation*}
\begin{diagram}
\chi:\ \ 0 & \rTo & L^\vee & \rTo & Q & \rTo & L & \rTo & 0
\end{diagram}
\end{equation*}
with 
\[
Q \cong \cok\left(\left(
   \begin{matrix} X & A \\ 0 & X^T\end{matrix}
   \right), \left(
   \begin{matrix} \adj(X) & -B_A \\ 0 & \adj(X)^T \end{matrix}
   \right)\right)
\]
for some alternating matrix $A$ over $S$ and its companion matrix
$B_A$.  We therefore have $\J_1(\chi) = I_1(A)K$ and $\J_2(\chi) =
I_1(B_A)S/I_{n-1}(X)$.  In particular, for each ideal of
$S/I_{n-1}(X)$ of the form $I_1(B_A)$, where $B_A$ is the companion
matrix for some alternating matrix $A$, there
exists an orientable MCM $R$-module $Q$ of rank $2$, and distinct
ideals yield nonisomorphic modules $Q$.  More precisely, we have the
following result.
\end{sit}

\begin{prop}\label{prop:reduce-by-2}
There is a surjective function from the isomorphism classes of
rank-two MCM $R$-modules appearing as the middle terms of extensions
of $L$ by $L^\vee$ to the set of principal ideals of the polynomial ring
in $(n-2)^2$ variables.
\end{prop}

\begin{proof}
Let $X'$ be the generic square matrix of size $n-2$, with entries
$x_{ij}'$, $1 \leq i,j \leq n-2$.  Let $S'=K[x'_{ij}]$ be the polynomial ring
over  $K$ in those indeterminates $x_{ij}'$, and define $\pi: S \to
S'$ by $\pi(x_{ij}) = x_{ij}'$ if $i,j \leq n-2$ and
$\pi(x_{ij})=0$ otherwise.  The $(n-1)$-minors of $X$ vanish
under $\pi$, so we obtain an induced epimorphism $\pi: S/I_{n-1}(X)
\to S'$.  Note that all $(n-2)$-minors of $X$ vanish under
$\pi$ as well, save $[n-1,n\;\widehat\mid\;n-1,n]$, which maps to
$\det X'$.  

Let $\chi \in \Ext_R^1(L,L^\vee)$.  Then $\chi$ is the image of an
alternating matrix $A$, and the ideal $\J_2(\chi) \subseteq
S/I_{n-1}(X)$ is generated by the entries of the companion matrix
$B_A$.  Again, $\J_2(\chi)$ depends only on the isomorphism class of
the middle term 
of $\chi$.  Recall (Theorem~\ref{thm:matrixmain}) that 
\[
b_{rs} = \sum_{k<l} a_{kl}(-1)^{r+s+k+l}[rs\;\widehat\mid\;kl]\,.
\]
The image of $\J_2(\chi)$ in $S'$, then, is generated by the single
element $\pi(a_{n-1,n})\cdot\det X'$.

Define $p: \Ext_R^1(L,L^\vee) \to \{\text{ideals of } S'\}$ by
$p(\chi) = (\pi(a_{n-1,n}))$.  Since $\det
X'$ is a nonzerodivisor in $S'$, $p(\chi)$ is a well-defined ideal
of $S'$.  Letting $A$ vary over all alternating matrices, we see
that $p$ is surjective, and by construction $p(\chi)$ depends only on
the isomorphism class of the middle term of $\chi$.
\end{proof}

\begin{cor} \label{cor:wild}
Let $X=(x_{ij})$ be the generic $(n \times n)$--matrix over the
field $K$, $n\geq 3$.  Let $R =K[x_{ij}]/(\det X)$ be the
generic determinantal hypersurface ring.  Then the rank-two
orientable MCM $R$-modules cannot be parametrized by the points of any
finite-dimensional algebraic variety over $K$.
\end{cor}

\begin{prob}
Our methods afford us no information in general about orientable
rank-two MCM $R$-modules generated by fewer than $n$ elements.  These
modules correspond to matrix factorizations $(\phi,\psi)$ of the
generic determinant of size $m<n$, with $\det \phi = (\det X)^2$ and
$\det \psi = (\det X)^{m-2}$ up to unit multiples.  When $n \leq 4$, no 
such module can exist: either $\det\phi$ or $\det \psi$ must be a unit
multiple of $\det X$, so must have cokernel among $\{L, L^\vee\}$, both
of which are $n$-generated. However, we do not know 
whether there exists a $4$-generated rank-two orientable MCM
$R$-module when $n=5$.  By the correspondence laid out in
\cite{Herzog-Kuhl:1989}, such a module would correspond to a
codimension-$3$ complete intersection ideal in $K[x_{1,1},
\dots,x_{5,5}]$ containing $\det X$ as a non-minimal  generator. 
\end{prob}

\section{Higher-order extensions}\label{section:higher-order}

In this final section we consider the higher-order extension theory of
the rank-one MCM modules over the generic determinant.  We shall see
that it is controlled by the ``half-trace'' of
Proposition~\ref{prop:matrixmain2} and Remark~\ref{rmk:half-trace}.

Maintain the notation of \ref{notation}. 
Recall from Corollary~\ref{cor:mainL-L} and
Proposition~\ref{prop:stable-ext-periodic} that we have natural
isomorphisms
\[
\Ext_R^2(L,L) \cong \uEnd_R(L) \cong S/I_{n-1}(X)\,.
\]
Dually, we also have
\[
\Ext_R^2(L^\vee,L^\vee) \cong \uEnd_R(L^\vee) \cong S/I_{n-1}(X)\,.
\]

\begin{theorem}\label{thm:Ext2-r}
Let $\chi \in \Ext_R^1(L,L^\vee)$ and $\chi' \in \Ext_R^1(L^\vee,L)$.  Let
$A = (a_{kl})$ and $A' = (a'_{kl})$ be alternating matrices
representing $\chi$ and $\chi'$, respectively, and let $r \in R$ and
$C$ be defined as in Theorem~\ref{thm:bifac}, so that
\[
r = -\sum_{k<l, u<v}(-1)^{u+v+k+l}a_{kl}
[uv\;\widehat\mid\;kl]a'_{uv}\ \in K\
\]
and we have
\[ 
A \adj(X) A' = r\cdot X^T + X^T C X^T\,.
\]
Then the image of $r$ in $S/I_{n-1}(X)$ represents both the Yoneda products
$\chi' \chi \in \Ext_R^2(L,L)$ and $-\chi\chi' \in \Ext_R^2(L^\vee,L^\vee)$.
\end{theorem}

\begin{sit} \label{sit:Yoneda}
Before beginning the proof, we recall the definition and computation
of the Yoneda product.  For extensions $\alpha \in \Ext_R^m(A,B)$ and
$\beta \in \Ext_R^n(B,C)$ over some ring $R$, represented by an
$m$-fold extension and an $n$-fold extension of $R$-modules,
respectively, the Yoneda product $\beta\alpha \in \Ext_R^{m+n}(A,C)$
is represented by the $(m+n)$-fold extension obtained by splicing the
representatives for $\alpha$ and $\beta$ together at their common
endpoint $B$.  This product is computed as follows
(\cite[p.91]{MacLane:homology}):
\begin{enumerate}
\item Choose preimages $\sigma \in \Hom_R(\syz^R_m(A), B)$ and $\rho \in
\Hom_R(\syz_n^R(B), C)$ for $\alpha$ and $\beta$;
\item Lift $\sigma$ to $\tilde\sigma \in
  \Hom_R(\syz^R_{m+n}(A),\syz^R_n(B))$\;
\item $\beta\alpha$ is the image of the composition $\rho\tilde\sigma$
  in $\Ext_R^{m+n}(A,C)$.
\end{enumerate}
\end{sit}

\begin{proof}[Proof of Theorem~\ref{thm:Ext2-r}]
To compute $\chi'\chi$, we first choose the natural preimages in
$\Hom_R(M,L^\vee)$ and $\Hom_R(M^\vee,L)$.  Specifically, $\cok(A,B_A) \in
\Hom_R(M,L^\vee)$ is a preimage for $\chi$, and $\cok(A',B_{A'}) \in
\Hom_R(M^\vee,L)$ is a preimage for $\chi'$.  We lift $\cok(A,B_A)$
naturally to $\cok(B_A,A) \in \Hom_R(L,M^\vee)$.  Then $\chi'\chi$ is
computed by the image in $\Ext_R^2(L,L)$ of $A'B_A \in
\Hom_R(\oF,\oF)$.  By Proposition~\ref{prop:matrixmain2}, we have
\[
B_A A' = r \cdot \id_n + CX^T\,.
\]
As $A', B_A$ are alternating, transposing both sides yields 
\[
A'B_A =  r \cdot \id_n + XC^T\,.
\]
The image of the term $XC^T$ in $\Ext_R^2(L,L)$ is zero, as $X C^T$
factors through $X$.  Thus $\chi'\chi$ is the image in $\Ext_R^2(L,L)$
of $r \cdot \id_n$, which agrees in $S/I_{n-1}(X)$ with $r$.

A symmetric calculation reveals that $r$ also represents $-\chi\chi'
\in \Ext_R^2(L^\vee,L^\vee)$.
\end{proof}

\begin{sit} 
The map 
\begin{alignat*}{2}
E = \Ext_R^1(L, L^\vee) \qquad \to & \qquad R_1 := S/I_{n-1}(X) \\
\chi = [A] \qquad \longmapsto & \qquad r(\tau(\chi),\chi) =
r(\tau(A),A) \\ 
& \qquad = -\sum_{k<l, 
  u<v}(-1)^{u+v+k+l}\tau(a_{kl})[uv\;\widehat\mid\;kl]a_{uv} 
\end{alignat*}
is $R_1$-quadratic.  We use this quadratic form to define the stable
$\Ext$-algebra. 
\end{sit}

\begin{defn}
The {\em stable $\Ext$-algebra of the rank-one MCM $R$-modules\/} is
the positively graded algebra $\calE$ with homogeneous components
\[
\calE^i = \uExt^i_R(L\oplus L^\vee,L\oplus L^\vee)\,,
\]
and multiplication induced by the Yoneda product.
\end{defn}

The graded components of $\calE$ depend only on parity, so we may
consider instead $\underline{\calE}:=\calE^0 \oplus \calE^1$ as a
graded algebra over $\bbbz/2\bbbz$.  The
structure of $\underline\calE$ can then be arranged as $(2\times
2)$--matrices:  
\begin{align*}
\calE^0 &= \uEnd_R(L\oplus L^\vee) \cong 
\left(\begin{matrix}R_1 & 0 \\ 0 & R_1 \end{matrix}\right)\,; \\
\calE^1 &= \uExt_R^1(L \oplus L^\vee, L\oplus L^\vee) \cong
\left(\begin{matrix}0 & E' \\ E & 0 \end{matrix}\right)
\end{align*}
with multiplication in $\calE^1$ defined by the quadratic form $r$.
Here we have observed that $\uHom_R(L,L^\vee)=0$ from the resolution
in \cite[Theorem 2.16]{Bruns-Vetter:1988}; see also the proof of
Theorem~\ref{thm:res-E}. 

Here is a summary of these observations.

\begin{theorem}\label{thm:ext-algebra}
With structure as defined above, the stable $\Ext$-algebra
$\underline\calE$ is a graded-commutative, $\bbbz/2\bbbz$-graded
algebra, with each homogeneous component an orientable MCM module of
rank two over $S/I_{n-1}(X)$.  Moreover, the multiplication yields
\[
\calE^0 / (\calE^1)^2 \cong 
\left(\begin{matrix}
S/I_{n-2}(X) & 0 \\ 0 & S/I_{n-2}(X)
\end{matrix}\right)\,,
\]
so that the quadratic form degenerates precisely on the singular locus
of $R_1$.
\end{theorem}


\begin{thebibliography}{10}

\bibitem{Avramov-Buchweitz}
L.L. Avramov and R.-O. Buchweitz, \emph{Support varieties and cohomology over
  complete intersections}, Invent. Math. \textbf{142} (2000), 285--318.

\bibitem{Baciu-Ene-Pfister-Popescu:2004}
C.~Baciu, V.~Ene, G.~Pfister, and D.~Popescu, \emph{Rank two {Cohen-Macaulay}
  modules over singularities of type $x_1^3+x_2^3+x_3^3+x_4^3$},
  J. Algebra \textbf{292} (2005), 447--491.

\bibitem{Bergman:2003}
G.M. Bergman, \emph{Can one factor the classical adjoint of a generic matrix?},
  Transformation Groups \textbf{11} (2006), 7--15.

\bibitem{Bruns:1975}
W.~Bruns, \emph{Die {Divisorenklassengruppe} der {R}estklassenringe von
  {P}olynomenringen nach {D}eterminantenidealen}, Revue {R}oumaine {M}ath.
  {P}ures {A}ppl. \textbf{20} (1985), 1109--1111.

\bibitem{BH}
W.~Bruns and J.~Herzog, \emph{{Cohen--Macaulay} rings}, Cambridge Stud. in Adv.
  Math., vol.~39, Cambridge University Press, Cambridge, 1993.

\bibitem{Bruns-Roemer-Wiebe:2005}
W.~Bruns, T.~{R\"omer}, and A.~Wiebe, \emph{Initial algebras of determinantal
  rings, {Cohen-Macaulay} and {U}lrich ideals}, Michigan Math. J. \textbf{53}
  (2005), 71--81.

\bibitem{Bruns-Vetter:1988}
W.~Bruns and U.~Vetter, \emph{Determinantal {Rings}}, Springer-Verlag, Berlin,
  1988, out of print; available at {\tt
  http://www.mathematik.uni-osnabrueck.de/preprints/shadow/calg9910.html}.

\bibitem{BGS}
R.-O. Buchweitz, G.-M. Greuel, and F.-O. Schreyer, \emph{{Cohen--Macaulay}
  modules on hypersurface singularities {II}}, Invent. Math. \textbf{88}
  (1987), 165--182.

\bibitem{DeConcini-Reichstein:2003}
C.~de~Concini and Z.~Reichstein, \emph{Nesting maps of {G}rassmannians},
  Rendiconti di {M}atematica {A}ccademia dei {L}incei \textbf{15} (2004),
  no.~s. 9, 109--118.

\bibitem{Eisenbud:1980}
D.~Eisenbud, \emph{Homological algebra on a complete intersection, with an
  application to group representations}, Trans. Amer. Math. Soc. \textbf{260}
  (1980), 35--64.

\bibitem{Eisenbud:book}
D.~Eisenbud, \emph{Commutative {Algebra} with a {View} {Toward} {Algebraic}
  {Geometry}}, Graduate {Texts} in {Mathematics}, vol. 150, Springer-Verlag,
  New York, 1995.

\bibitem{M2}
D.~Grayson and M.~Stillman, \emph{Macaulay 2, a software system for research in
  algebraic geometry}, Available at {\tt http://www.math.uiuc.edu/Macaulay2/}.

\bibitem{Gulliksen-Negard:1972}
T.H. Gulliksen and O.G. Neg{\aa}rd, \emph{Un complexe r\'esolvant pour certains
  id\'eaux d\'eterminantiels}, C. R. Acad. Sci. Paris S\'er. A-B \textbf{274}
  (1972), A16--A18.

\bibitem{Herzog:1978}
J.~Herzog, \emph{Ringe mit nur endlich vielen {Isomorphieklassen} von maximalen
  unzerlegbaren {Cohen--Macaulay Moduln}}, Math. Ann. \textbf{233} (1978),
  21--34.

\bibitem{Herzog-Kuhl:1989}
J.~Herzog and M.~K{\"u}hl, \emph{Maximal {C}ohen-{M}acaulay modules over
  {G}orenstein rings and {B}ourbaki-sequences}, Commutative Algebra and
  Combinatorics (Kyoto, 1985), Adv. Stud. Pure Math., vol.~11, North-Holland,
  Amsterdam, 1987, pp.~65--92.

\bibitem{Herzog-Trung:1992}
J.~Herzog and N.V. Trung, \emph{Gr\"obner bases and multiplicity of
  determinantal and {P}faffian ideals}, Adv. Math. \textbf{96} (1992), no.~1,
  1--37.

\bibitem{Hochster-Eagon:1971}
M.~Hochster and J.~Eagon, \emph{{Cohen--Macaulay} rings, invariant theory, and
  the generic perfection of determinantal loci}, Amer. J. Math. \textbf{53}
  (1971), 1020--1058.

\bibitem{Ile:2004}
R.~Ile, \emph{Deformation theory of rank one maximal {Cohen--Macaulay} modules
  on hypersurface singularities and the {Scandinavian} complex}, Compositio
  Math. \textbf{140} (2004), 435--446.

\bibitem{MacLane:homology}
S.~Mac Lane, \emph{Homology}, Classics in Mathematics, Springer-Verlag, Berlin,
  1995, reprint of the 1975 edition.

\bibitem{Miyata}
T.~Miyata, \emph{Note on direct summands of modules}, J. Math. Kyoto Univ.
  \textbf{7} (1967), 65--69.

\end{thebibliography}

\providecommand{\bysame}{\leavevmode\hbox to3em{\hrulefill}\thinspace}
\providecommand{\MR}{\relax\ifhmode\unskip\space\fi MR }
\providecommand{\MRhref}[2]{%
  \href{http://www.ams.org/mathscinet-getitem?mr=#1}{#2}
}
\providecommand{\href}[2]{#2}

\end{document}